\documentclass{article}
\usepackage{amsmath,amsfonts,amssymb,amsthm}
\usepackage{graphicx}
\usepackage{subfig}
\usepackage[title]{appendix}
\usepackage{fancyhdr}
\usepackage{listings}
\usepackage{color, colortbl}
\usepackage{xcolor}
\usepackage{multirow}
\usepackage{geometry}
\usepackage{array,multirow}
\usepackage{mathtools}
\usepackage[]{algorithm2e}
\usepackage{mathrsfs}
\usepackage{enumitem}
\usepackage{booktabs}
\usepackage{adjustbox}
\usepackage{authblk}
\usepackage{url}
\usepackage{hyperref}

\newtheorem{theorem}{Theorem}[section]
\newtheorem{definition}[theorem]{Definition}

\newtheorem{assump}[]{}

\newcounter{algo}[section]
\renewcommand{\thealgo}{\thesection.\arabic{algo}}
\newcommand{\algo}[3]{\refstepcounter{algo}
\begin{center}
\framebox[\textwidth]{
\parbox{0.95\textwidth} {\vspace{\topsep}
{\bf Algorithm \thealgo : #2}\label{#1}\\
\vspace*{-\topsep} \mbox{ }\\
{#3} \vspace{\topsep} }}
\end{center}}

\theoremstyle{plain}
\title{Convergent least-squares optimisation methods for variational data assimilation}
\author[4]{C. Cartis}
\author[1,*]{M. H. Kaouri}
\author[1,2,3]{A. S. Lawless}
\author[1,2,3]{N. K. Nichols}

\affil[1]{\footnotesize Department of Mathematics and Statistics, University of Reading, UK}
\affil[2]{\footnotesize Department of Meteorology, University of Reading, UK}
\affil[3]{\footnotesize National Centre for Earth Observation, Reading, UK}
\affil[4]{\footnotesize Mathematical Institute, University of Oxford, UK}
\affil[*]{\footnotesize Corresponding author.
Address: Department of Mathematics and Statistics, University of Reading, PO Box 220, Reading, RG6 6AX, UK. Email: mahakaouri@gmail.com}

\begin{document}

\maketitle
\begin{abstract}
Data assimilation combines prior (or background) information with observations to estimate the initial state of a dynamical system over a given time-window. A common application is in numerical weather prediction where a previous forecast and atmospheric observations are used to obtain the initial conditions for a numerical weather forecast. In four-dimensional variational data assimilation (4D-Var), the problem is formulated as a nonlinear least-squares problem, usually solved using a variant of the classical Gauss-Newton (GN) method. However, we show that GN may not converge if poorly initialised. In particular, we show that this may occur when there is greater uncertainty in the background information compared to the observations, or when a long time-window is used in 4D-Var allowing more observations. The difficulties GN encounters may lead to inaccurate initial state conditions for subsequent forecasts. To overcome this, we apply two convergent GN variants (line search and regularisation) to the long time-window 4D-Var problem and investigate the cases where they locate a more accurate estimate compared to GN within a given budget of computational time and cost. We show that these methods are able to improve the estimate of the initial state, which may lead to a more accurate forecast.
\end{abstract}

\newline
\textbf{Highlights:}
\begin{itemize}
\item Poor initialisation of Gauss-Newton method may result in failure to converge.
\item Safeguarded Gauss-Newton improves initial state estimate within limited time/cost.
\item Results using twin experiments with long time-window and chaotic Lorenz models.
\item Apply state of the art least-squares convergence theory to data assimilation.
\item Improvements to initial state estimate may lead to a more accurate forecast.
\end{itemize}

\section{Introduction}

Four-dimensional variational data assimilation (4D-Var) aims to solve a nonlinear least-squares problem that minimizes the error in a prior estimate of the initial state of a dynamical system together with the errors between observations and numerical model estimates of the states of the system over time. In Numerical Weather Prediction (NWP), 4D-Var is used to estimate the initial conditions for a weather forecast \cite{le1986variational}. The 4D-Var scheme is able to incorporate information from a previous forecast along with observations over both temporal and spatial domains, weighted by the uncertainties in the prior and the observations. From a Bayesian point of view the solution is the \textit{maximum a posteriori} estimate of the initial state \cite{nichols2010mathematical}. The nonlinear least-squares objective function is minimized using an iterative method. The quality of the estimate and the subsequent forecast depends on how accurately the 4D-Var problem is solved within the time and computational cost available. \\
 
In this paper, we investigate the application of \textit{globally convergent} optimisation methods to the 4D-Var problem; such methods use safeguards to guarantee convergence from an arbitrary initial estimate by ensuring a sufficient, monotonic/strict decrease in the objective function at each iteration. We focus on the strong-constraint 4D-Var problem where we assume that the numerical model of the system perfectly represents the true dynamics of the system or the model errors are small enough to be neglected. This results in the formulation of variational data assimilation as an unconstrained nonlinear least-squares problem and is employed by many operational meteorological centres \cite{rabier2005overview}, including the Meteorological Service of Canada \cite{gauthier2007extension}, the European Centre for Medium-range Weather Forecasting (ECMWF) \cite{courtier1994strategy, rabier1998extended} and the Met Office \cite{rawlins2007met}. \\

Ideally in large-scale unconstrained optimisation, we seek a fast rate of convergence, which can be achieved in nondegenerate cases using a Newton-type method. However, these methods require the use of second order derivatives of the objective function, which are too costly to compute and store operationally. Therefore, optimisation methods that approximate the high order terms, such as limited memory Quasi-Newton \cite{gilbert1989some, liu1989limited, shanno1980remark, zou1993numerical}, Inexact Newton \cite{dembo1982inexact}, Truncated Newton \cite{le2002second, wang1995truncated}, Adjoint Newton \cite{wang1998adjoint}, Hessian-free Newton \cite{daescu2003analysis}, Gauss-Newton \cite{dennis1996numerical, nocedal2006numerical} and Approximate Gauss-Newton \cite{gratton2007approximate} methods have been considered. To compute efficiently the first derivatives of the objective function required by these techniques, the adjoint of the numerical model is generally used \cite{le1986variational}. More recently, optimisation methods that do not require the first derivatives of the objective function are being examined to avoid the development and maintenance costs associated with using the adjoint \cite{gratton2014derivative}. Alternative data assimilation techniques that use ensemble methods to approximate the objective function gradients, rather than using the adjoint, are also being investigated \cite{bannister2017review, liu2008ensemble}. \\

The incremental 4D-Var technique, used commonly in operational centres, approximately solves a sequence of linear least-squares problems and has been shown to be equivalent to the Gauss-Newton (GN) method under standard conditions \cite{lawless2005investigation}. In the GN (or incremental) method the linearized problem is solved in an inner loop of the algorithm; the solution to the nonlinear problem is then updated in an outer loop and the problem is re-linearized. The accuracy with which the inner loop is solved is known to affect the convergence of the outer loop \cite{laroche1998validation, lawless2005investigation}. In our work, we focus on the convergence of the outer loop, where the exact gradient is used (as is the case when an adjoint is available) and we assume that the inner loop linear least-squares problem is solved exactly. Furthermore, we use a variable transformation usually applied in operational 4D-Var to precondition the optimisation problem, see \cite{bannister2008review2}. \\

A general drawback of the GN method is that given a poor initialisation, it is not guaranteed to converge to a solution, known as the ‘analysis’ state, of the 4D-Var problem \cite{dennis1996numerical}. In NWP, the initial guess for the minimisation is generally chosen to be the predicted initial state from a previous forecast, known as the ‘prior’ or ‘background’ state. However, for some applications of 4D-Var this choice may not be a good enough estimate of the analysis. We show that in such cases, the GN minimisation procedure may fail to converge. There are three main strategies that safeguard GN and make it convergent from an arbitrary initial guess: line-search, regularisation and trust-region \cite{conn2000trust, nocedal2006numerical}. These can all be regarded as variants of the original Levenberg-Marquardt algorithm for solving nonlinear least-squares problems \cite{levenberg1944method, marquardt1963algorithm}. In this work we investigate GN method with regularisation (REG) and compare its performance to GN with backtracking Armijo line-search (LS) and GN alone, applied to the preconditioned 4D-Var problem when there is limited computational time and evaluations available, such as in NWP. \\

In previous work, the use of a line-search strategy in combination with a Quasi-Newton approach was implemented in the ECMWF NWP system to solve the 4D-Var problem and was found to improve the minimisation of the objective function \cite{gilbert1989some, rabier1998extended}. This method uses the Wolfe line-search conditions \cite{wolfe1969convergence} to safeguard the convergence. The Wolfe conditions require the use of additional evaluations of the objective function and its gradient, however, which is computationally costly. Here we instead use the Armijo condition \cite{armijo1966minimization}, which only requires additional evaluations of the objective function and not the gradient. We pair GN with backtracking Armijo line-search and use a fixed number of computational evaluations to guarantee a reduction in the outer loop objective function (assuming the inner loop is solved to a high accuracy). We compare this method to the GN method and to the GN method safeguarded by quadratic regularisation (REG), using a simple, inexpensive updating strategy. \\

Using two test models within the 4D-Var framework, we show that where there is more uncertainty in the background information compared to the observations, the GN method may fail to converge, yet the convergent methods, LS and REG, are able to improve the estimate of the analysis. Assimilation over long time windows is of particular interest. We use accuracy profiles to show numerically that in the long time-window case and in cases where there is higher uncertainty in the background information versus the observations, the globally convergent methods are able to solve more problems than GN in the limited cost available. By `solve' we mean satisfying a criterion requiring a reduction in the objective function within a set number of evaluations. We also show the effect that poor background information has on the quality of the estimate obtained. We consider the case where the background information is highly inaccurate compared to the observations and find that the convergence of all three methods is improved when more observations are included along the time-window. Finally, for the case where GN performs well, we recommend further research into the parameter updating strategies used within the globally convergent methods. \\

The structure of this paper is organised as follows. In Section \ref{vardasection} we outline the strong-constraint 4D-Var problem as a nonlinear least-squares problem and the GN method that is frequently used to solve it. In Section \ref{globconvsec} we outline the globally convergent methods used within this paper. In Section \ref{expdessec} we describe the experimental design including the dynamical models used. In Section \ref{numexp} we present the numerical results obtained when applying GN and the globally convergent methods to the 4D-Var problem with different features. Finally, we conclude our findings in Section \ref{concsec}. In an appendix we detail the proofs of convergence for the REG and LS methods.

\section{Variational data assimilation} 
\label{vardasection}

\subsection{4D-Var: least-squares formulation}
In four-dimensional variational data assimilation (4D-Var), the analysis $\mathbf{x}_0^a \in \mathbb{R}^n$ is obtained by minimising a objective function consisting of two terms: the background term and the observation term, namely;

\begin{equation}
\label{4dvar}
\mathcal{J}(\mathbf{x}_0) = \frac{1}{2}(\mathbf{x}_0-\mathbf{x}_0^b)^T\mathbf{B}_0^{-1}(\mathbf{x}_0-\mathbf{x}_0^b)+\frac{1}{2}\sum_{i=0}^N (\mathbf{y}_i-\mathcal{H}_i(\mathbf{x}_i))^T\mathbf{R}_i^{-1}(\mathbf{y}_i-\mathcal{H}_i(\mathbf{x}_i)).
\end{equation}

The background term measures the difference between the initial state of the system and the background state vector $\mathbf{x}_0^b \in \mathbb{R}^n$, which contains prior information. The observation term measures the difference between information from observations at times $t_i$ in the observation vector $\mathbf{y}_i \in \mathbb{R}^{p_i}$ and the model state vector $\mathbf{x}_i \in \mathbb{R}^{n}$ at the same time through use of the observation operator $\mathcal{H}_i: \mathbb{R}^n \rightarrow \mathbb{R}^{p_i}$ that maps from the model state space to the observation space. Both terms are weighted by their corresponding covariance matrices to represent the uncertainty in the respective measures, the background error covariance matrix $\mathbf{B} \in \mathbb{R}^{n \times n}$ and the observation error covariance matrices at times $t_i$, $\mathbf{R}_i \in \mathbb{R}^{p_i \times p_i}$, which are assumed to be symmetric positive definite. We note that observations are distributed both in time and space and there are usually fewer observations available than there are state variables so $p < n$, where $p = \sum_{i=0}^N p_i$. The 4D-Var objective function \eqref{4dvar} is subject to the nonlinear dynamical model equations which contain the physics of the system
\begin{equation}
\label{4dmodelcons}
\mathbf{x}_{i} = \mathcal{M}_{0,i}(\mathbf{x}_0),
\end{equation}
where the nonlinear model $\mathcal{M}_{0,i}: \mathbb{R}^n \rightarrow \mathbb{R}^n$ evolves the state vector from the initial time point $t_0$ to the time point $t_i$.\\ 

We precondition the 4D-Var problem using a variable transform, which has been shown to improve the conditioning of the variational optimisation problem \cite{haben2011conditioning, haben2011conditioninginc}. To be able to use the negative square root of $\mathbf{B}$ in our variable transformation, we first require the assumption that the matrix $\mathbf{B}$ is full rank. This assumption is satisfied for our choices of $\mathbf{B}$ in Section \ref{numexp}. We define a new variable $\mathbf{v}$ to be,
\begin{equation}
\label{CVTvar}
    \mathbf{v} = \mathbf{B}^{-1/2}(\mathbf{x}_0-\mathbf{x}_0^b).
\end{equation}
The 4D-Var objective function can then be written in terms of $\mathbf{v}$, known as the control variable in data assimilation (DA), and minimised with respect to this instead. Furthermore, by including the model information within the objective function, we are able to write the constrained optimisation problem \eqref{4dvar}-\eqref{4dmodelcons} in the form of an unconstrained optimisation problem and apply the minimisation methods described later in this paper. The preconditioned 4D-Var objective function is given by
\begin{equation}
\label{4dvarp}
\mathcal{J}(\mathbf{v}) = \frac{1}{2}\mathbf{v}^T\mathbf{v}+\frac{1}{2}\sum_{i=0}^N (\mathbf{y}_i-\mathcal{H}_i(\mathcal{M}_{0,i}(\mathbf{B}^{1/2}\mathbf{v} +\mathbf{x}_0^b)))^T\mathbf{R}_i^{-1}(\mathbf{y}_i-\mathcal{H}_i(\mathcal{M}_{0,i}(\mathbf{B}^{1/2}\mathbf{v} +\mathbf{x}_0^b))).
\end{equation}
We note that the function \eqref{4dvarp} is continuously differentiable if the operators $\mathcal{H}_i$ and $\mathcal{M}_{0,i}$ are continuously differentiable. To save both computational cost and time in 4D-Var, tangent linear approximations of the nonlinear operators in \eqref{4dvarp} are use in the inner loop \cite{courtier1994strategy}. The tangent linear model (TLM) and tangent linear observation operator are usually derived by linearising the discrete nonlinear model equations.\\

In a nonlinear least-squares problem, the function $\mathcal{J}: \mathbb{R}^n \rightarrow \mathbb{R}$ has a special form, as defined in the following,
\begin{equation}
\label{resform}
\min_{\mathbf{v}} \mathcal{J}(\mathbf{v}) = \frac{1}{2} \|\mathbf{r}(\mathbf{v})\|_2^2,
\end{equation}
where $\mathbf{r}(\mathbf{v}) = [r_1(\mathbf{v}),...,r_{n+p}(\mathbf{v})]^T$ and each $r_j:\mathbb{R}^n \rightarrow \mathbb{R}$, for $j=1, 2, \dots, n+p$, is referred to as a residual. In \eqref{resform}, $\| . \|_2$ denotes the $l_2$-norm, which will be used throughout this paper. Equation \eqref{4dvarp} is equivalent to \eqref{resform} where the residual vector $\mathbf{r}(\mathbf{v}) \in \mathbb{R}^{(n+p)}$ and its Jacobian $\mathbf{J}(\mathbf{v})$ are given by
\begin{equation}
\label{precondresJac}
    \mathbf{r}(\mathbf{v}) = \begin{pmatrix}
    \mathbf{v}\\
    \mathbf{R}_0^{-1/2}(\mathbf{y}_0-\mathcal{H}_0(\mathbf{B}^{1/2}\mathbf{v} +\mathbf{x}_0^b))\\
    \mathbf{R}_1^{-1/2}(\mathbf{y}_1-\mathcal{H}_1(\mathcal{M}_{0,1}(\mathbf{B}^{1/2}\mathbf{v} +\mathbf{x}_0^b)))\\
    \vdots \\
    \mathbf{R}_N^{-1/2}(\mathbf{y}_N-\mathcal{H}_N(\mathcal{M}_{0,N}(\mathbf{B}^{1/2}\mathbf{v} +\mathbf{x}_0^b))
    )
    \end{pmatrix} \text{ and } 
    \mathbf{J}(\mathbf{v}) = \begin{pmatrix}
    \mathbf{I}\\
    -\mathbf{R}_0^{-1/2}\mathbf{H}_0\mathbf{B}^{1/2}\\
    -\mathbf{R}_1^{-1/2}\mathbf{H}_1\mathbf{M}_{0,1}\mathbf{B}^{1/2}\\
    \vdots\\
    -\mathbf{R}_N^{-1/2}\mathbf{H}_N\mathbf{M}_{0,N}\mathbf{B}^{1/2}\\
    \end{pmatrix},
\end{equation}

where 
\begin{equation}
    \mathbf{M}_{0,i} = \frac{\partial\mathcal{M}_{0,i}}{\partial\mathbf{v}}\big|_{\mathcal{M}_{0,i}(\mathbf{B}^{1/2}\mathbf{v} +\mathbf{x}_0^b)} \text{ and } \mathbf{H}_i = \frac{\partial\mathcal{H}_0}{\partial\mathbf{v}}\big|_{\mathcal{M}_{0,i}(\mathbf{B}^{1/2}\mathbf{v} +\mathbf{x}_0^b)}
\end{equation}
are the Jacobian matrices of the model operator $\mathcal{M}_{0,i}$ and observation operator $\mathcal{H}_i$ respectively, where $\mathbf{M}_{0,i} \in \mathbb{R}^{n \times n}$ is the tangent linear of $\mathcal{M}_{0,i}$ and $\mathbf{H}_i \in \mathbb{R}^{p_i \times n}$ is the tangent linear of $\mathcal{H}_i$ \cite{nichols2010mathematical}. In practice, an adjoint method is used to calculate the gradient of \eqref{4dvarp}, defined as
\begin{equation}
\label{4dvgrad}
\nabla \mathcal{J}(\mathbf{v}) = \mathbf{J}(\mathbf{v})^T\mathbf{r}(\mathbf{v}).
\end{equation}
The Hessian is the matrix of second-order partial derivatives of \eqref{4dvarp},
\begin{equation}
    \label{4dvarHess}
    \nabla^2 \mathcal{J}(\mathbf{v}) = \mathbf{J}(\mathbf{v})^T\mathbf{J}(\mathbf{v}) + \sum\limits_{j=1}^{n+p} r_j(\mathbf{v}) \nabla^2 r_j(\mathbf{v}).
\end{equation}
In data assimilation, the second-order terms in \eqref{4dvarHess} are often difficult to calculate in the time and cost available and too large to store, and so one cannot easily use Newton-type methods for 4D-Var. Therefore, a first-order approximation to the Hessian of the objective function \eqref{4dvarp} is used, resulting in a GN method, and is given by
\begin{equation}
\label{precondhess}
    \mathbf{S} = \mathbf{J}(\mathbf{v})^T\mathbf{J}(\mathbf{v}) = \mathbf{I} + \sum_{i=0}^N \mathbf{B}^{1/2}\mathbf{M}_{0,i}^T\mathbf{H}_i^T\mathbf{R}_i^{-1}\mathbf{H}_i\mathbf{M}_{0,i}\mathbf{B}^{1/2}, 
\end{equation}
which is, by construction, full rank and symmetric positive definite. The condition number in the $l_2$-norm of \eqref{precondhess}, $\kappa(\mathbf{S})$, is the ratio of its largest and smallest eigenvalues and is related to the number of iterations used for the linear minimisation problems in 4D-Var and how sensitive the estimate of the initial state is to perturbations of the data. We can use $\kappa(\mathbf{S})$ to indicate how quickly and accurately the optimisation problem can be solved \cite{golub2012matrix}. 

\subsection{4D-Var implementation}
The incremental 4D-Var method, which was first proposed for practical implementation of the NWP problem in \cite{courtier1994strategy}, has been shown to be equivalent to the GN method when an exact TLM is used in the inner loop. When an approximate TLM is used, the method is equivalent to an inexact GN method \cite{gratton2007approximate, lawless2005investigation}. A summary of the GN method is given next.

\vskip 3pt
\algo{GNalg}{GN algorithm applied to \eqref{resform} \cite{dennis1996numerical}.}{
\textbf{Step \boldmath$0$: Initialisation}. Given $\mathbf{v}^{(0)} \in \mathbb{R}^n$ and some stopping criteria. Set $k=0$. 
\vskip 2pt
\textbf{Step \boldmath$1$: Check stopping criteria.} While the stopping criteria are not satisfied, do:
\vskip 2pt
\indent\hspace{1cm} \textbf{Step \boldmath$2$: Step computation.} Compute a step $\mathbf{s}^{(k)}$ that satisfies 
    \begin{equation}
    \label{GNsk}
        \mathbf{J}(\mathbf{v}^{(k)})^T\mathbf{J}(\mathbf{v}^{(k)})\mathbf{s}^{(k)} = -\mathbf{J}(\mathbf{v}^{(k)})^T\mathbf{r}(\mathbf{v}^{(k)}).
    \end{equation}
\vskip 2pt
\indent\hspace{1cm} \textbf{Step \boldmath$3$: Iterate update.} Set $\mathbf{v}^{(k+1)} = \mathbf{v}^{(k)} + \mathbf{s}^{(k)}$, $k := k+1$ and go to Step 1.
\vskip 2pt
}
\vskip 2pt

In Algorithm \ref{GNalg}, the updated control variable $\mathbf{v}^{(k+1)}$ is computed by finding a step $\mathbf{s}^{(k)}$ that satisfies \eqref{GNsk}, which is known as the preconditioned linearised subproblem. By substituting $\mathbf{v}^{(k+1)}$ into \eqref{CVTvar} and rearranging, we obtain the current estimate $\mathbf{x}_0^{(k+1)}$ of the initial state to the original nonlinear 4D-Var problem.\\ 

To reduce the computational cost in large DA systems and to solve the DA problem in real time, the series of problems \eqref{GNsk} can be solved approximately in the inner loop using iterative optimisation methods such as Conjugate Gradient (CG) where a limited number of CG iterations are allowed and an exact or approximate $\mathbf{J}$ is used \cite{gratton2007approximate}. We do not focus on this here and assume that \eqref{GNsk} is solved exactly. \\ 

We note that the step calculation \eqref{GNsk} uniquely defines $\mathbf{s}^{(k)}$, and $\mathbf{s}^{(k)}$ is a descent direction when $\mathbf{J}(\mathbf{v})$ is full column rank. This is the case in 4D-Var as the Jacobian, $\mathbf{J}(\mathbf{v})$ in \eqref{precondresJac} is full column rank due to the presence of the identity matrix, thus ensuring that $\mathbf{s}^{(k)}$ is a descent direction. \\

The definitions of two solution types, namely, local and global minima, are stated in Appendix \ref{convthmappend}, along with a brief explanation of the local convergence property of GN. Although the GN method benefits from local convergence properties, convergence can only be guaranteed if the initial guess $\mathbf{v}^{(0)}$ of the algorithm is in some neighbourhood around an (unknown) local solution $\mathbf{v}^*$, that is, convergence from an arbitrary initial guess is not guaranteed \cite{dennis1996numerical}. Even if the GN method does converge, it may not necessarily converge to the global minimum due to the fact that multiple local minima of a nonlinear least-squares objective function may exist.\\ 

GN has no way of adjusting the length of the step $\mathbf{s}^{(k)}$ and hence, may take steps that are too long and fail to decrease the objective function value and thus to converge, see Example 10.2.5 in \cite{dennis1996numerical} and later in Section \ref{numexp} where the poor performance of GN is demonstrated. As GN only guarantees local convergence, we are interested in investigating methods that converge when $\mathbf{v}^{(0)}$ is far away from a local minimiser $\mathbf{v}^{*}$. We refer to these methods as `globally convergent'. Mathematical theory on global strategies can be found in \cite{nocedal2006numerical} and \cite{dennis1996numerical}. Two globally convergent methods are GN with line search and GN with quadratic regularisation, which use a strategy within the GN framework to achieve convergence to a stationary point given an arbitrary initial guess by adjusting the length of the step. These methods will be presented in the next section. 

\section{Globally convergent methods}
\label{globconvsec}
Within this section, we outline the two globally convergent algorithms that we apply in Section \ref{numexp} to the preconditioned 4D-Var problem.

\subsection{Gauss-Newton with line search (LS)}
A line search method aims to restrict the step $\mathbf{s}^{(k)}$ in \eqref{GNsk} so as to guarantee a decrease in the value of $\mathcal{J}$. Within our work, an inexact line search method known as the backtracking-Armijo (bArmijo) algorithm is used within the inner loop of GN to find a step length $\alpha>0$ that satisfies the Armijo condition \cite{armijo1966minimization}. The Gauss-Newton with backtracking-Armijo line search (LS) method is as follows.

\vskip 3pt
\algo{LSalg}{LS algorithm applied to \eqref{resform} \cite{nocedal2006numerical}.}{
\textbf{Step \boldmath$0$: Initialisation}. Given $\mathbf{v}^{(0)} \in \mathbb{R}^n$, $\tau \in (0,1)$ and $\beta \in (0,1)$ and $\alpha_0 > 0$ and some stopping criteria. Set $k=0$. 
\vskip 2pt
\textbf{Step \boldmath$1$: Check stopping criteria.} While the stopping criteria are not satisfied, do:
\vskip 2pt
\indent\hspace{1cm} \textbf{Step \boldmath$2$: Step computation.} Compute a step $\mathbf{s}^{(k)}$ that satisfies 
    \begin{equation}
    \label{LSsk}
        \mathbf{J}(\mathbf{v}^{(k)})^T\mathbf{J}(\mathbf{v}^{(k)})\mathbf{s}^{(k)} = -\mathbf{J}(\mathbf{v}^{(k)})^T\mathbf{r}(\mathbf{v}^{(k)})
    \end{equation}
\indent\hspace{1cm} and set $\alpha^{(k)} = \alpha_0$.
\vskip 2pt
\indent\hspace{1cm} \textbf{Step \boldmath$3$: Check Armijo condition.} While the following (Armijo) condition is not satisfied
    \begin{equation}
    \label{bArmijoalg}
        \mathcal{J}(\mathbf{v}^{(k)} + \alpha^{(k)}\mathbf{s}^{(k)}) \leq \mathcal{J}(\mathbf{v}^{(k)}) + \beta\alpha^{(k)}{\mathbf{s}^{(k)}}^T \nabla \mathcal{J}(\mathbf{v}^{(k)}),
    \end{equation}
\indent\hspace{1cm} do: 
\vskip 2pt
\indent\hspace{2cm} \textbf{Step \boldmath$4$: Shrink stepsize.} Set $\alpha^{(k)} := \tau\alpha^{(k)}$ and go to Step 3. 
\vskip 2pt
\indent\hspace{1cm} \textbf{Step \boldmath$5$: Iterate update.} Set $\mathbf{v}^{(k+1)} = \mathbf{v}^{(k)} + \alpha^{(k)}\mathbf{s}^{(k)}$, $k := k+1$ and go to Step 1.
\vskip 2pt
}
\vskip 2pt

In Algorithm \ref{LSalg}, the control parameter $\beta$ in \eqref{bArmijoalg} is typically chosen to be small (see \cite{nocedal2006numerical}). The step equation \eqref{LSsk} is the same as the GN step equation \eqref{GNsk}; thus when $\alpha^{(k)}=1$, the GN and LS iterates coincide at (the same) point $\mathbf{v}^{(k)}$. The use of condition \eqref{bArmijoalg} in this method ensures that the accepted steps produce a sequence of strictly decreasing function values given $\triangledown \mathcal{J}(\mathbf{v}^{(k)})^T \mathbf{s}^{(k)}<0$. This latter condition is satisfied by $\mathbf{s}^{(k)}$ defined in \eqref{LSsk} whenever $\mathcal{J}(\mathbf{v}^{(k)})$ is full column rank (which is the case here) as mentioned in Section \ref{vardasection} \cite{nocedal2006numerical}.\\

Despite its global convergence property (see Appendix \ref{LSappendix}), the LS method has some disadvantages. We remark that the use of the step length $\alpha^{(k)}$ may sometimes unnecessarily shorten the step $\mathbf{s}^{(k)}$, slowing down the convergence. Furthermore, LS may be computationally costly due to the need to calculate the value of the function $\mathcal{J}$ each time $\alpha^{(k)}$ is adjusted, although more sophisticated updating strategies for $\alpha$ may be used to try to  reduce this effect.\\

Other line search strategies are possible such as Wolfe, Goldstein-Armijo and more \cite{nocedal2006numerical}, but they are more involved and potentially more computationally costly. As LS requires the re-evaluation of the outer loop objective function each time it adjusts its line search parameter, its applicability to real systems has been in doubt due to the computational cost limitations in 4D-Var \cite{rabier1998extended}. In Section \ref{numexp}, we show that given the same cost as the GN method, the LS method can in some cases, better minimise the preconditioned 4D-Var objective function.\\

\subsection{Gauss-Newton with regularisation (REG)}
The GN method may also be equipped with a globalisation strategy by including a regularisation term $\gamma^{(k)}\mathbf{s}^{(k)}$ in the step calculation \eqref{GNsk} of Algorithm \ref{GNalg}. This ensures that the accepted steps produce a sequence of monotonically decreasing function values. This is a common variation of the GN method known as the Levenberg-Marquardt method, proposed in \cite{levenberg1944method} and \cite{marquardt1963algorithm}. The effect of the regularisation parameter $\gamma^{(k)}$ is to implicitly control the length of the step $\mathbf{s}^{(k)}$. Increasing $\gamma^{(k)}$ shortens the steps, thus increasing the possibility that the procedure will decrease the objective function in the next iteration. The REG method can be summarised as follows.

\vskip 3pt
\algo{REGalg}{REG algorithm applied to \eqref{resform} \cite{more1978levenberg}.}{
\textbf{Step \boldmath$0$: Initialisation}. Given $\mathbf{x}^{(0)} \in \mathbb{R}^n$, $1 > \eta_2 \geq \eta_1 > 0$, $\gamma^{(0)} > 0$ and some stopping criteria. Set $k=0$. 
\vskip 2pt
\textbf{Step \boldmath$1$: Check stopping criteria.} While the stopping criteria are not satisfied, do:
\vskip 2pt
\indent\hspace{1cm} \textbf{Step \boldmath$2$: Step computation.} Compute a step $\mathbf{s}^{(k)}$ that satisfies
    \begin{equation}
    \label{regsk}
        \left( \mathbf{J}(\mathbf{v}^{(k)})^T\mathbf{J}(\mathbf{v}^{(k)}) +  \gamma^{(k)}\mathbf{I}\right)\mathbf{s}^{(k)} = -\mathbf{J}(\mathbf{v}^{(k)})^T\mathbf{r}(\mathbf{v}^{(k)}).
    \end{equation}
\vskip 2pt
\indent\hspace{1cm} \textbf{Step \boldmath$3$: Iterate update.} 
Compute the ratio
    \begin{equation}
    \label{rhoalg}
        \rho^{(k)} = \frac{\mathcal{J}(\mathbf{v}^{(k)}) - \mathcal{J}(\mathbf{v}^{(k)} + \mathbf{s}^{(k)})}{\mathcal{J}(\mathbf{v}^{(k)}) - m(\mathbf{s}^{(k)})},
    \end{equation}
\indent\hspace{1cm} where 
    \begin{equation}
\label{mapprox}
m(\mathbf{s}^{(k)}) = \frac{1}{2}\|\mathbf{J}(\mathbf{v}^{(k)})\mathbf{s}^{(k)} + \mathbf{r}(\mathbf{v}^{(k)})\|_2^2 + \frac{1}{2}\gamma^{(k)}\|\mathbf{s}^{(k)}\|_2^2.
\end{equation}
\indent\hspace{1cm} Set 
    \begin{equation}
        \mathbf{v}^{(k+1)} = \begin{cases}
    \mathbf{v}^{(k)} + \mathbf{s}^{(k)}, & \text{if $\rho^{(k)} \geq \eta_1$}\\
    \mathbf{v}^{(k)}, & \text{otherwise}.
  \end{cases}
    \end{equation}
\vskip 2pt
\indent\hspace{1cm} \textbf{Step \boldmath$4$: Regularisation  parameters update.} Set 
    \begin{equation}
    \label{gammaupdate}
        \gamma^{(k+1)} = \begin{cases}
    \frac{1}{2}\gamma^{(k)}, & \text{if $\rho^{(k)} \geq \eta_2$} \text{   (very successful iteration)} \\
    \gamma^{(k)}, & \text{if $\eta_1 \leq \rho^{(k)} < \eta_2$} \text{   (successful iteration)} \\
    2\gamma^{(k)}, & \text{otherwise,} \text{   (unsuccessful iteration)}
  \end{cases}
    \end{equation} 
\indent\hspace{1cm} Let $k := k+1$ and go to Step 1.
\vskip 2pt
}
\vskip 2pt

As in Algorithms \ref{GNalg} and \ref{LSalg}, the step equation \eqref{regsk} is solved directly in the numerical experiments in Section \ref{numexp}. We note that when $\gamma^{(k)} =0$ in \eqref{regsk}, the REG step in \eqref{regsk} is the same as the GN step in \eqref{GNsk}. By comparing \eqref{regsk} with \eqref{GNsk}, we are able to see how the REG step differs from the GN step. The diagonal entries of the Hessian of the 4D-Var objective function \eqref{4dvarp} are increased by the regularisation parameter $\gamma^{(k)}$ at each iteration of the REG method. The method is able to vary its step between a GN and a gradient descent step by adjusting $\gamma^{(k)}$ (see \cite{nocedal2006numerical}) but may be costly due to the need to calculate the value of the function $\mathcal{J}$ to assess the step. Note that other choices of the factors $\frac{1}{2}$ and $2$ for updating $\gamma^{(k)}$ in \eqref{gammaupdate} are possible and even more sophisticated variants for choosing $\gamma^{(k)}$ have been proposed. The proof of global convergence of the REG method is presented in Appendix \ref{REGappendix}. 

\section{Experimental design}
\label{expdessec}

Before evaluating the GN, LS and REG methods numerically, we first explain the experimental design.\\ 

Twin experiments are commonly used to test DA methods. They use error statistics that satisfy the DA assumptions as well as synthetic observations generated by running the nonlinear model forward in time to produce a reference state (not generally a local minimum of \eqref{resform}). Within this section, we define our choices for the twin experimental design. We begin by briefly outlining two commonly used dynamical models, which are sensitive to initial conditions (chaotic nature), a property shared with NWP models.

\subsection{Models}
\textbf{Lorenz 1963 model (L63)} Proposed in \cite{lorenz1963deterministic}, the Lorenz 63 model (L63) is a popular experimental dynamical system that represents meteorological processes using a simple model. The model consists of three nonlinear, ordinary differential equations given as
\begin{equation}
\label{lorenz63}
\begin{aligned}
\frac{dx}{dt} &= \sigma(y - x), \\
\frac{dy}{dt} &= x(\rho - z) - y, \\
\frac{dz}{dt} &= xy - \beta z,
\end{aligned}
\end{equation}
where the state vector consists of $n=3$ time-dependent variables $\mathbf{x} = [x(t), y(t), z(t)]^T \in \mathbb{R}^3$. The scalar parameters are chosen to be $\sigma = 10$, $\rho = \frac{8}{3}$ and $\beta = 28$, making the system chaotic. A second-order Runge-Kutta method is used to discretise the model equations using a time step $\Delta t = 0.025$. \\
\newline
\textbf{Lorenz 1996 model (L96)} Another popular experimental system is the atmospheric Lorenz 96 model (L96) \cite{lorenz1996predictability} given by the following $n$ equations, 
\begin{equation}
    \label{lorenz96}
    \frac{dx_j}{dt} = -x_{j-2} x_{j-1} + x_{j-1} x_{j+1} - x_{j} + F,
\end{equation}
where $j=1,2,\hdots,n$ is a spatial coordinate. For a forcing term $F=8$ and $n=40$ state variables, the system is chaotic \cite{lorenz1996predictability}. The variables are evenly distributed over a circle of latitude of the Earth with $n$ points with a cyclic domain and a single time unit is equivalent to approximately 5 atmospheric days. A fourth-order Runge-Kutta method is used to discretise the model equations using a time step $\Delta t = 0.025$ (approximately 3 hours). \\

For both the L63 and L96 models, the time-window length $t_a$ is varied in the numerical experiments in Section \ref{efftwl}. We will now outline how we formulate the twin experiments, beginning with generating the reference state.

\subsection{Twin experiments}

The reference state at time $t_0$, $\mathbf{x}^{ref}_0$ is used as the basis of a twin experiment in the definition of the background state (the initial guess for the optimisation algorithms) as well as to generate the observations using a nonlinear model run called the `nature' run. We begin by explaining how we obtain $\mathbf{x}^{ref}_0$.\\
\newline 
\textbf{Reference state} A vector of length $n$ is drawn from the uniform distribution and used as the initial vector of state variables $\mathbf{x}^{rand}$. For the L63 model, $\mathbf{x}^{rand}$ is integrated forward using a second-order Runge-Kutta method, which is spun-up over 1000 time steps to obtain the reference state on the model attractor for the L63 twin experiments, $\mathbf{x}^{ref}_0 \in \mathbb{R}^{3}$. This is the same for the L96 model except a fourth-order Runge-Kutta method is used to obtain $\mathbf{x}^{ref}_0 \in \mathbb{R}^{40}$. The reference state at time $t_0$, $\mathbf{x}^{ref}_0$ can then be used to obtain the full nonlinear model trajectory. \\

We next explain how we obtain the background state vector used within our twin experiments to be used as the initial guess for the optimisation algorithms.\\
\newline
\textbf{Background} In 4D-Var, the initial guess for the optimisation algorithm is taken to be the background state at time $t_0$, $\mathbf{x}_0^b$, which incorporates information from previous forecasts. In our experiments, the background state vector $\mathbf{x}^{b}_0$ is generated by adding Gaussian noise
\begin{equation}
\label{bgerr}
    \mathbf{\varepsilon_b} \sim \mathcal{N}(0, \mathbf{B}),
\end{equation}
to the reference state at time $t_0$, $\mathbf{x}^{ref}_0$. For the background error covariance matrix, we choose $\mathbf{B} = \sigma_b^2 \mathbf{I}_n$ where $\sigma_b^2$ is the background error variance. The standard deviations of the errors from the reference state for each model are based on the average order of magnitude of the entries of $\mathbf{x}^{ref}_0$. For the L63 experiments, $\sigma_b^2 = 0.25, 1, 6.25$ and $25$ represent a $5\%, 10\%, 25\%$ and $50\%$ error respectively. Similarly for the L96 experiments we set $\sigma_b^2 = 0.0625, 0.25, 1.5625$ and $6.25$. \\

As previously mentioned, we generate synthetic observations from a nonlinear model run using the reference state at time $t_0$, $\mathbf{x}^{ref}_0$. We next describe the choices we made when specifying these observations. \\
\newline
\textbf{Observations} We consider both the spatial and temporal locations of the observations. We assume that for both models observations of single state variables are taken and $\mathbf{H}_i$ are the exact observation operators at times $t_i$ used to map to observation space. For the L63 model, we consider where we have $p=2$ observations, one of $x$ and one of $z$ per observation location in time. For the L96 model, we consider where we have an observation of the first half of the state variables per observation location in time. This choice mimics what we may expect in reality where we have more observations concentrated in one part of the globe. For both models, we first consider where there is only one set of observations at time $N$ (Nobs1) and then show the effect of using more observations along the time-window in later experiments. We use imperfect observations where the observations $\mathbf{y}_i$ are generated by adding Gaussian noise 
\begin{equation}
\label{obserr}
    \mathbf{\varepsilon_o} \sim \mathcal{N}(0, \mathbf{R}_i),
\end{equation}
to $\mathbf{H}_i \mathbf{x}^{ref}_i$ for each observation location in time. For the observation error covariance matrix we choose $\mathbf{R}_i = \sigma_o^2 \mathbf{I}_{p}$ where $\sigma_o^2$ is the observation error variance. We expect the problem \eqref{4dvarp} to be more ill-conditioned, thus difficult to solve accurately, when the ratio 
\begin{equation}
    \label{sigmaratio}
    \frac{\sigma_b}{\sigma_o}
\end{equation}
is large \cite{haben2011conditioning, haben2011conditioninginc}. The ratio \eqref{sigmaratio} controls the influence of the observation term in the preconditioned objective function \eqref{4dvarp}. For all experiments, we set the standard deviation of the observation error to be $10\%$ of the average order of magnitude of the entries of $\mathcal{H}(\mathbf{x}^{ref}_i)$ for both models. For the L63 model, this is $\sigma_o^2 = 1$ and for the L96 model, this is $\sigma_o^2 = 0.25$. We vary the background error variance $\sigma_b^2$ above and below $\sigma_o^2$ such that the ratio \eqref{sigmaratio} varies. This can be thought of as having more confidence in the observations compared to background when $\sigma_b > \sigma_o$ and vice versa. Furthermore, as the initial guess is set to be the background state vector, which is dependent on the value of $\sigma_b$, by varying $\sigma_b^2$ we are essentially varying the initial guess of the algorithms, thus eliminating starting point bias from our results \cite{beiranvand2017best}. It is important to recall here that under certain conditions, the GN method is known for its fast convergence properties when in close vicinity to a local minimum, see \cite{dennis1996numerical}. By choosing a small value of $\sigma_b^2$, we expect the performance of GN to beat that of both LS and REG as it does not require the adjustment of the additional parameters $\alpha^{(k)}$ and $\gamma^{(k)}$. Also, when assuming that the observations are more accurate than the background, the use of more observation locations in time means that we are constraining the estimate of the initial state more tightly to the reference state in the twin experiment design. The effect this has on the convergence of the optimisation methods will be investigated. We next outline the algorithmic choices we have made.

\subsection{Algorithmic choices}

\textbf{Stopping criteria} We now outline the criteria used to terminate Algorithms \ref{GNalg}, \ref{LSalg} and \ref{REGalg}. Due to the limited time and computational cost available in practice, the GN method is not necessarily run to convergence and a stopping criterion is used to limit the number of iterations. Each calculation of the residual vector $\mathbf{r}(\mathbf{v})$ requires the non-linear model to be run forward to obtain the state at each observation location in time. This can then be used to calculate the value of the objective function. Furthermore, one run of the adjoint model is required to calculate the gradient.\\ 

To reduce computational cost in practical implementations of 4D-Var, the inner loop problem is solved at a lower resolution than the outer loop problem \cite{fisher2009data}. However, as the dimension of the problems used within this paper are relatively small compared to DA systems in practice, we solve the full resolution inner loop problem using the full resolution residual and Jacobian given in \eqref{precondresJac} and solve the inner loop problem using MATLAB's backslash operator where an appropriate solver is chosen according to the properties of the Hessian matrix $\nabla^2 \mathcal{J}(\mathbf{v})$ (see \cite{MATLABmldivide2021} for more details). The limit on the total number of function and Jacobian evaluations is achieved by using the following criterion
\begin{equation}
    \label{funcjaclim}
    k_J+l \leq \tau_e, 
\end{equation}
where $k_J$ is the total number of Jacobian evaluations (which is equivalent to the number of outer iterations $k$ in 4D-Var), $l$ is the total number of function evaluations and $\tau_e$ is the tolerance. The tolerance $\tau_e$ can be chosen according to the maximum number of evaluations desired. We note that for GN, $k_J = l$ as the method requires as many Jacobian evaluations as function evaluations. However, for both LS and REG there could be more than one function evaluation per Jacobian evaluation since for unsuccessful steps, the Jacobian is not updated so $k_J \leq l$.\\

To ensure that the algorithms are stopped before the function values stagnate, the following common termination criterion based on the relative change in the function at each iteration is also used
\begin{equation}
    \label{relfunccrit}
    \frac{|\mathcal{J}(\mathbf{v}^{(k-1)}) - \mathcal{J}(\mathbf{v}^{(k)})|}{1 +\mathcal{J}(\mathbf{v}^{(k)})} \leq \tau_s,
\end{equation}
for $k \geq 1$, where $\tau_s$ is the tolerance, chosen to be $10^{-5}$. The criterion \eqref{relfunccrit} is used throughout Section \ref{numexp} unless indicated otherwise.\\ 

We expect the norm of the gradient of the objective function, $\|\nabla \mathcal{J}(\mathbf{v}^{(k)})\|$ to be close to zero at a stationary point. The following termination criterion will be used in Section \ref{behdivGN} to identify whether or not a given method has located a stationary point
\begin{equation}
    \label{normgradcrit}
    \|\nabla \mathcal{J}(\mathbf{v}^{(k)})\| \leq 10^{-5}.
\end{equation}
\newline
\textbf{Parameter choices}
For the LS method, we choose $\alpha_0 = 1$ so that the first step assessed by the bArmijo rule is the GN step. We set $\beta = 0.1$ and to adjust the step length, $\tau = 0.5$.\\ 

For the REG method, we select the initial regularisation parameter to be $\gamma^{(0)} = 1$ so that the condition in Algorithm \ref{REGalg}, $\gamma^{(0)}>0$, is satisfied and the REG step differs from the GN step. Furthermore, we choose $\eta_1 = 0.1$ and $\eta_2 = 0.9$ to assess how well the model \eqref{mapprox} approximates the true function value at the next iteration.\\ 

For all three optimisation methods, we set $\tau_e = 8, 100$ or $1000$ depending on the experiment. The choice of $\tau_e = 8$ comes from that which is used operationally in the ECMWF Integrated Forecasting System \cite{ECMWFnewsletter158}, whereas the choice of $\tau_e = 100$ or $1000$ is used to measure the performance of the optimisation methods when closer to convergence.\\

In order to best present our results, we use accuracy profiling described as follows. \\
\newline
\textbf{Accuracy profiles}
An \textit{accuracy profile} shows the proportion of problems a given method can solve within a fixed amount of work ($\tau_e$) and a given tolerance ($\tau_f$) of the change in the function value \cite{more2009benchmarking}. To ensure the robustness of our results, we apply the three optimisation methods to a series of $n_r$ randomly generated problems, where the randomness occurs through the background and observation error vectors, $\mathbf{\varepsilon_b}$ and $\mathbf{\varepsilon_o}$. For each realisation, a new $\mathbf{\varepsilon_b}$ and $\mathbf{\varepsilon_o}$ are generated from their respective distributions, \eqref{bgerr} and \eqref{obserr}. The following criterion proposed in \cite{more2009benchmarking} is used to flag that an estimate of the initial state has been obtained by an optimisation method
\begin{equation}
    \label{funcxastop}
    \frac{\mathcal{J}(\mathbf{v}_0^{(l)}) - \mathcal{J}(\mathbf{v}_0^t)}{\mathcal{J}(\mathbf{v}_0^{(0)}) - \mathcal{J}(\mathbf{v}_0^t)} \leq \tau_f, 
\end{equation}
where $\mathbf{v}_0^t$ is a solution of \eqref{4dvarp} referred to as the `truth' and $\tau_f$ is the tolerance. The measure \eqref{funcxastop} compares the optimality gap $\mathcal{J}(\mathbf{v}_0^{(l)}) - \mathcal{J}(\mathbf{v}_0^{t})$ relative to the best reduction $\mathcal{J}(\mathbf{v}_0^{(0)}) - \mathcal{J}(\mathbf{v}_0^t)$ \cite{more2009benchmarking}. This ensures that the 4D-Var problem is only flagged as solved by the optimisation method once the value of the objective function is within some error ($\tau_f$) of the truth.\\ 
For our problems, the truth is unknown. We only know that, due to the nonlinearity of the 4D-Var problem, there may exist many values of $\mathbf{v}_0$ that could minimise \eqref{4dvarp} locally. We are interested in the estimate $\mathbf{v}_0^t$ that gives the greatest reduction in \eqref{4dvarp} that any of the three methods can obtain. Therefore, we set the truth to be the $\mathbf{v}_0^{(l)}$ obtained by any of the three methods that gives the smallest function value within the given number of evaluations. Using this criterion allows us to benchmark the methods against each other using accuracy profiles.\\ 

For each experiment, we plot the proportion of the same $n_r=100$ realisations solved by each method against the relative accuracy obtained, $\tau_f$. The relative accuracy obtained is varied using $\tau_f = 10^{-i}$, where $i = 0, 0.01, 0.02, \dots, 5$.

\section{Numerical results}
\label{numexp}
In this section, we present the results when applying GN, LS and REG using the experimental design described in the previous section. We begin by conducting experiments showing the effect of the length of the assimilation time-window on the convergence of the three methods. 

\subsection{Effect of time-window length}
\label{efftwl}

\begin{figure}[h!]
\centering
 \subfloat[][L63, $t_a = 0.05$.]{\includegraphics[width=0.25\textwidth]{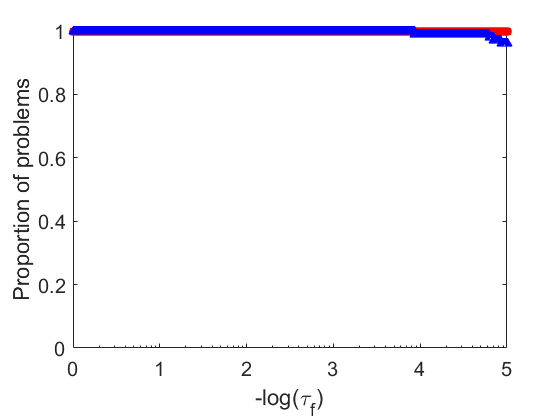}}
 \subfloat[][L63, $t_a = 0.1$.]{\includegraphics[width=0.25\textwidth]{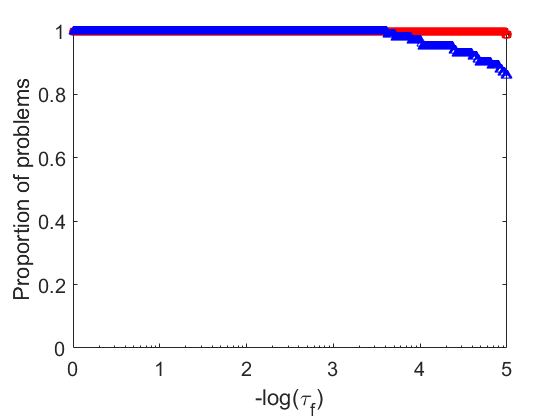}}
 \subfloat[][L63, $t_a = 0.2$.]{\includegraphics[width=0.25\textwidth]{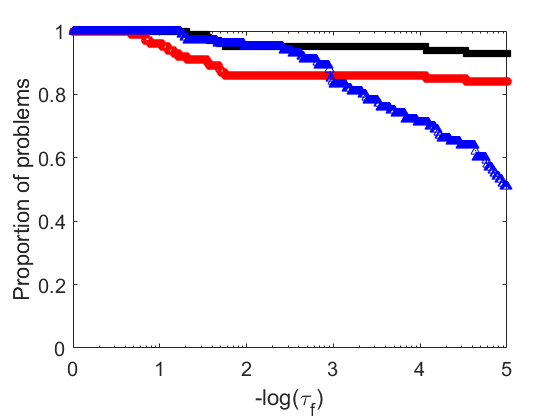}}
 \subfloat[][L63, $t_a = 1$.]{\includegraphics[width=0.25\textwidth]{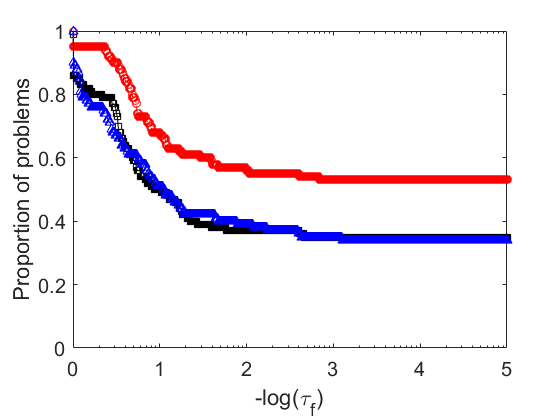}}\\
 
  \subfloat[][L96, $t_a = 0.05$.]{\includegraphics[width=0.25\textwidth]{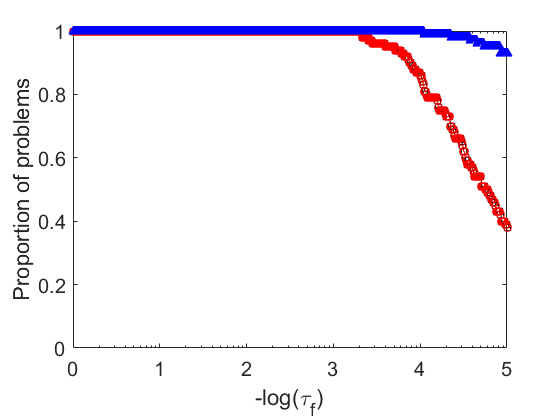}}
 \subfloat[][L96, $t_a = 0.1$.]{\includegraphics[width=0.25\textwidth]{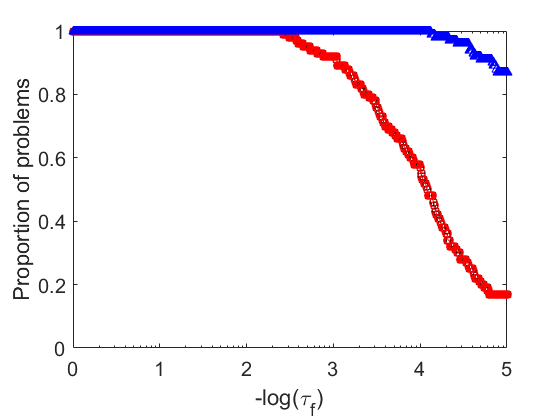}}
 \subfloat[][L96, $t_a = 0.2$.]{\includegraphics[width=0.25\textwidth]{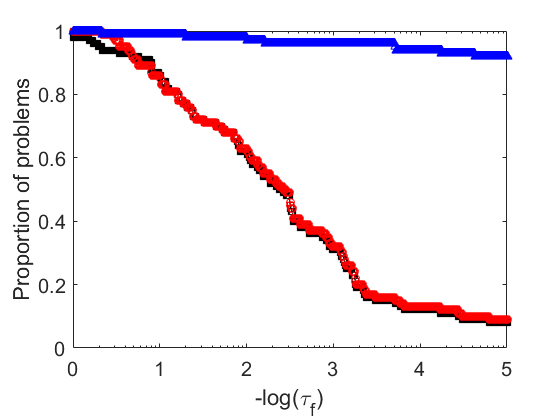}}
 \subfloat[][L96, $t_a = 1$.]{\includegraphics[width=0.25\textwidth]{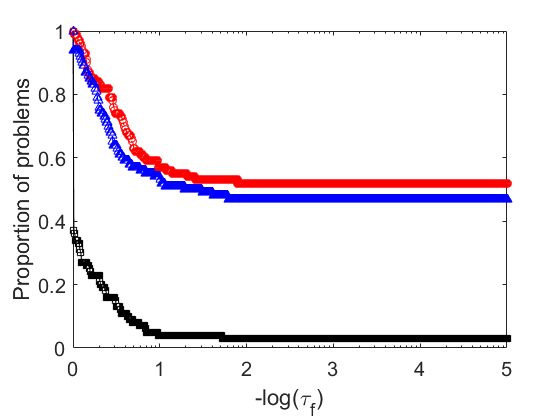}}
 \caption{Accuracy profiles for the GN (black), LS (red) and REG (blue) methods applied to the L63 and L96 problems using different time-window lengths $t_a$. These show the proportion of $n_r=100$ problems solved by each of the methods against the specified accuracy $-\log(\tau_f)$ when $\tau_e = 8$. The GN line is below the LS line in (a), (b), (e), (f) and (g).}
 \label{fig:timewindowplots}
\end{figure}

We produce accuracy profiles for different time-window lengths to understand the effect this has on the convergence of each method while limiting the number of function and Jacobian evaluations to $\tau_e = 8$. We choose a background error of $50\%$ and an observation error of $10\%$ so that the ratio \eqref{sigmaratio} is large relative to the other cases we consider. For both the L63 and L96 models, we consider both short and long time-window lengths of 6 hours ($t_a = 0.05$), 12 hours ($t_a = 0.1$), 1 day ($t_a = 0.2$) and 5 days ($t_a = 1$) with the results shown in Figure \ref{fig:timewindowplots}.\\

From Figure \ref{fig:timewindowplots}, we see that as the length of the time-window of both the L63 and L96 problems is increased, the performance of the GN, LS and REG methods suffers.\\

For the L63 problems, Figures \ref{fig:timewindowplots}(a) and \ref{fig:timewindowplots}(b) show that GN and LS perform similarly and solve more problems to the highest accuracy than REG. However, as $\tau_f$ is increased, REG is solving all of the problems, so the REG estimate must be close to that of GN and LS. In Figure \ref{fig:timewindowplots}(c), both LS and REG solve fewer problems compared to GN, even for relatively large choices of $\tau_f$. However, there is a choice of $\tau_f$ where all three methods are solving all problems, again indicating that the LS and REG estimates are close to the GN estimate. The initial guess for the three methods (the background) appears to be close enough to the solution and so the GN step is able to attain a sufficient decrease in the objective function as predicted by its local convergence properties. LS and REG are inadvertently shortening the GN step, which is a good step in the short time-window case. As we know, LS and REG need to adjust their respective parameters, $\alpha^{(k)}$ and $\gamma^{(k)}$ to attain GN's fast local convergence, so LS and REG are requiring more evaluations than GN to achieve the same result. For the L96 short time-window results in Figures \ref{fig:timewindowplots}(e), \ref{fig:timewindowplots}(f) and \ref{fig:timewindowplots}(g), this is not the case. In fact, REG is outperforming GN and LS and it appears that LS is mimicking the behaviour of GN quite closely as the GN step is attaining a sufficient decrease in the objective function. However the decrease that the REG step is achieving appears to be much greater for the L96 problems. Therefore, REG is able to solve a greater number of problems within a higher level of accuracy, which explains the difference between the L63 results in Figures \ref{fig:timewindowplots}(a), \ref{fig:timewindowplots}(b) and the L96 results in \ref{fig:timewindowplots}(e) and \ref{fig:timewindowplots}(f).\\

The long time-window results for the L63 and L96 problems are shown in Figures \ref{fig:timewindowplots}(d) and  \ref{fig:timewindowplots}(h), respectively. In both figures, LS is outperforming GN. For the L63 problems, the performance of GN does not differ much from the performance of REG. However, comparing the performance of GN in \ref{fig:timewindowplots}(c) with \ref{fig:timewindowplots}(d), we can see that performance of GN has deteriorated greatly when increasing the length of the time-window. In fact, in the results where even longer time-windows are used (not included here), LS and REG outperform the GN method for the L63 problems, as in \ref{fig:timewindowplots}(h). \\

For the remainder of our experiments, we set $t_a = 1$ in order to consider a long time-window case only, as this is where we expect to see the greatest benefit from the globally convergent methods.

\subsection{Behaviour of methods and stagnation of GN}
\label{behdivGN}
In order to gain an understanding of how the globally convergent methods, LS and REG, contend with GN, we next demonstrate the behaviour of GN, LS and REG when applied to typical preconditioned 4D-Var L63 and L96 problems, where the background error is large and the time-window length is long. \\

\begin{figure}[h!]
\centering
 \subfloat[][L63, Nobs1, $\sigma_b^2 = 25$, $\mathbf{B} = \sigma_b^2\mathbf{I}$, $t_a = 1$, $\tau_e = 100$.]{\includegraphics[width=0.5\textwidth]{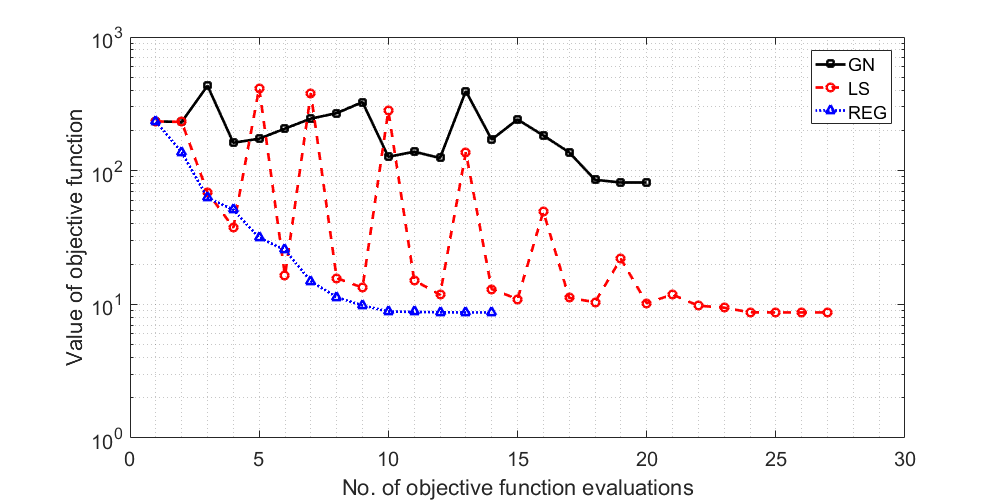}}
 \subfloat[][L96, Nobs1, $\sigma_b^2 = 6.25$, $\mathbf{B} = \sigma_b^2\mathbf{I}$, $t_a = 1$ $\tau_e = 100$.]{\includegraphics[width=0.5\textwidth]{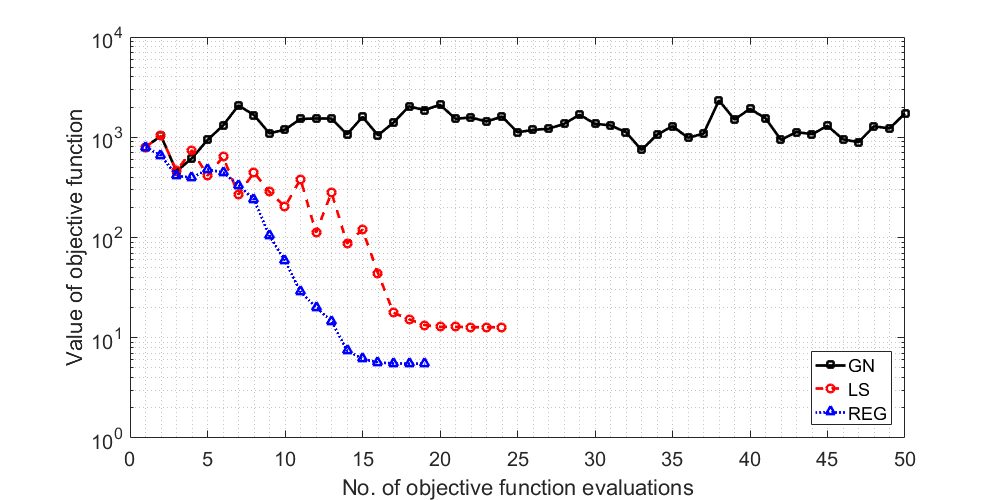}}\\
 \caption{Convergence plots showing the value of the objective function at each iteration (including unsuccessful iterations) of the GN (black), LS (red) and REG (blue) methods when applied to a L63 problem (a) and a L96 problem (b).}
 \label{fig:divplots}
\end{figure}

Figure \ref{fig:divplots} shows the convergence plots for two typical realisations when using the GN, LS and REG methods to obtain a solution to the preconditioned 4D-Var problem with the L63 and L96 models. In this figure, the total number of function and Jacobian evaluations allowed is set to $\tau_e = 100$ for both the L63 and the L96 problems to see if any progress is made beyond the number of evaluations allowed in practice. We recall that GN updates the gradient \eqref{4dvgrad} when the function \eqref{4dvarp} is updated, so there are as many function evaluations as Jacobian evaluations. However, both LS and REG only update the Jacobian on successful iterations when there is a reduction in the objective function. Therefore, the total number of evaluations used by each of the methods could consist of a different combination of function and Jacobian evaluations. As in Section \ref{efftwl}, we set the ratio \eqref{sigmaratio} to be large. It is in this case that we are able to best demonstrate the benefit of the globally convergent methods, LS and REG. In Figure \ref{fig:divplots}, we set $\tau_s = 10^{-3}$ to ensure that the methods stop before the function values stagnate. As Figure \ref{fig:divplots} includes function evaluations for both successful and unsuccessful step calculations, it is natural to see jumps in the function values of LS and REG while their parameters, $\alpha^{(k)}$ and $\gamma^{(k)}$ are being adjusted to guarantee a reduction in the function.\\

For the L63 problems (Figure \ref{fig:divplots}(a)), all three methods stop when the relative change in the function criterion \eqref{relfunccrit} is satisfied and before the limit on the total number of function and Jacobian evaluations \eqref{funcjaclim} is met. Table \ref{table:L63} supports this figure by showing the algorithmic output for each of the GN, LS and REG methods when two different stopping criteria are used. From these results, we see that both LS and REG stop at the same function value, although REG requires fewer evaluations to do so, and that GN is converging towards a larger value of the objective function \eqref{4dvarp} than LS and REG. By instead stopping on the criterion \eqref{normgradcrit} and setting $\tau_e = 1000$, we see in Table \ref{table:L63} that all three methods are still making progress on the gradient and iterate level, indicating that the methods are in fact locating stationary points despite a small change in the function value beyond those shown in Figure \ref{fig:divplots}. \\

\begin{table}[!ht]
\small
\caption{Table of algorithmic output when applying, GN, LS and REG to a typical realisation of the L63 problems, corresponding to Figure \ref{fig:divplots}(a).}
\centering
\begin{tabular}{ c c c c c c c}
\hline\hline
 Criteria & Method & $l$ & $k_J$ & $\mathcal{J}(\mathbf{v}^{(k_J)})$ & $\| \mathbf{v}^{(k_J)} - \mathbf{v}^{(k_J-1)} \|$ & $\| \nabla \mathcal{J}(\mathbf{v}^{(k_J)})\|$ \\ [0.5ex]
\hline
& GN & 20 & 20 & 81.55 & 0.42 & 86.35 \\ 
\eqref{relfunccrit} & LS & 27 & 14 & 8.69 & 0.03 & 5.18 \\ 
& REG & 14 & 14 & 8.69 & 0.05 & 1.00 \\ 
\hline
& GN & 101 & 101 & 78.87 & $3.54^{-8}$ & $8.47^{-6}$ \\ 
\eqref{normgradcrit} & LS & 43 & 27 & 8.69 & $8.21^{-7}$ & $8.31^{-6}$ \\ 
& REG & 66 & 66 & 8.69 & $7.34^{-7}$ & $9.24^{-6}$ \\ 
\hline
\hline
\end{tabular}
\label{table:L63}
\end{table}

For the L96 problems (Figure \ref{fig:divplots}(b)), LS and REG stop when \eqref{relfunccrit} is satisfied and before \eqref{funcjaclim} is satisfied, whereas GN only satisfies \eqref{funcjaclim}. Table \ref{table:L96} supports this figure by showing the algorithmic output for each of the GN, LS and REG methods when two different stopping criteria are used. From these results, we see that both GN and LS are stopping at a larger value of the objective function \eqref{4dvarp} than REG. Recall that the norm of the gradient criterion \eqref{normgradcrit} can be used to identify whether or not a given method has located a stationary point. The values of $\|\nabla \mathcal{J}(\mathbf{v}^{(k_J)})\|$ for LS and REG when the relative change in the function criterion \eqref{relfunccrit} is used are much smaller than that of GN. However, when we instead use the norm of the gradient criterion \eqref{normgradcrit} and limit the number of iterations to $\tau_e = 1000$, the methods stop on the limit of the number of iterations. Therefore, our results do not indicate that the estimates of LS and REG may indeed be stationary points of the objective function as they did for the L63 problems. However, LS and REG are are able to make some improvement (REG more so than LS) on the gradient norm level, unlike GN, which appears to fluctuate at gradient level, even after $\tau_e = 1000$ evaluations. \\ 

\begin{table}[!ht]
\small
\caption{Table of algorithmic output when applying, GN, LS and REG to a typical realisation of the L96 problems, corresponding to Figure \ref{fig:divplots}(b).}
\centering
\begin{tabular}{ c c c c c c c}
\hline\hline
 Criteria & Method & $l$ & $k_J$ & $\mathcal{J}(\mathbf{x}^{(k_J)})$ & $\| \mathbf{v}^{(k_J)} - \mathbf{v}^{(k_J-1)} \|$ & $\| \nabla \mathcal{J}(\mathbf{v}^{(k_J)})\|$ \\ [0.5ex]
\hline
& GN & 50 & 50 & 1728.99 & 20.02 & 5758.47\\ 
\eqref{relfunccrit} & LS & 24 & 14 & 12.72 & 0.07 & 10.09\\  
& REG & 19 & 16 & 5.52 & 0.08 & 1.89 \\  
\hline
& GN & 500 & 500 & 960.32 & 15.88 & 8015.13 \\ 
\eqref{normgradcrit} & LS & 967 & 32 & 12.71 & 0 & 10.09 \\ 
& REG & 967 & 32 & 5.51 & 0 & 0.03 \\ 
\hline
\hline
\end{tabular}
\label{table:L96}
\end{table}

Table \ref{table:L96} shows that as LS and REG iterate beyond what is shown in Figure \ref{fig:divplots}(b), there is very little change in the value of the cost function, despite making some change on the iterate and/or gradient level. The effect of rounding error means that although we see progress made, the function value may remain stagnant because of limitations in computer precision and because of the conditioning of the problem. The condition number of the Hessian $\kappa(\mathbf{S})$ can be used to indicate the accuracy we could be able to achieve. In our work, both the L63 and L96 problems are well-conditioned. \\

The observed behaviour in this section is partly due to the fact that there is no mechanism in GN to force it to converge as there is in LS and REG. The benefit of these mechanisms is clearly shown in Figure \ref{fig:divplots}(b) where the GN method is stagnating while the LS and REG methods are converging, further motivating our investigation of these methods.

\subsection{Effect of background error variance}
\label{effbgerr}

\begin{figure}[h!]
\centering
 \subfloat[][L63, $50\%$]{\includegraphics[width=0.25\textwidth]{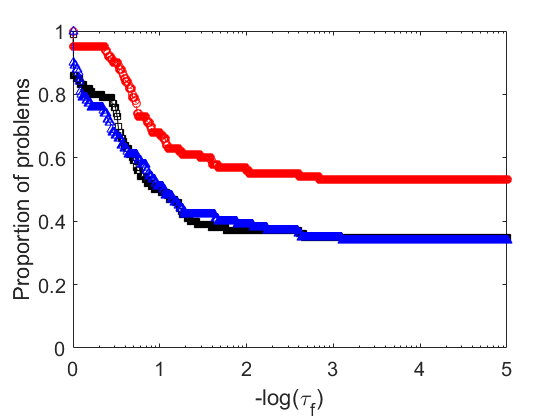}}
 \subfloat[][L63, $25\%$]{\includegraphics[width=0.25\textwidth]{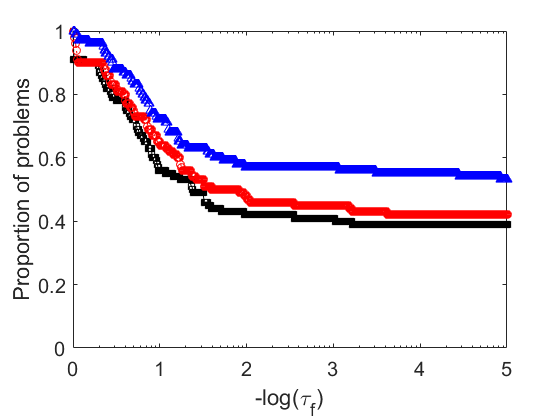}}
 \subfloat[][L63, $10\%$]{\includegraphics[width=0.25\textwidth]{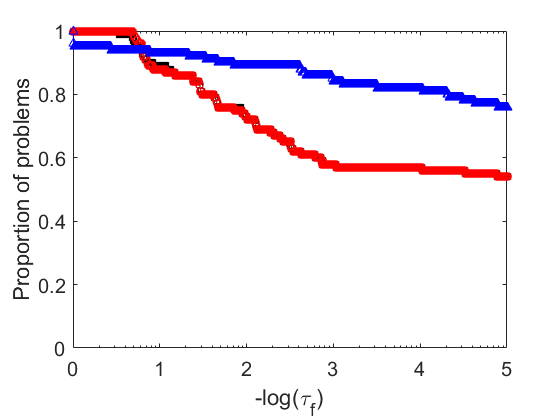}}
 \subfloat[][L63, $5\%$]{\includegraphics[width=0.25\textwidth]{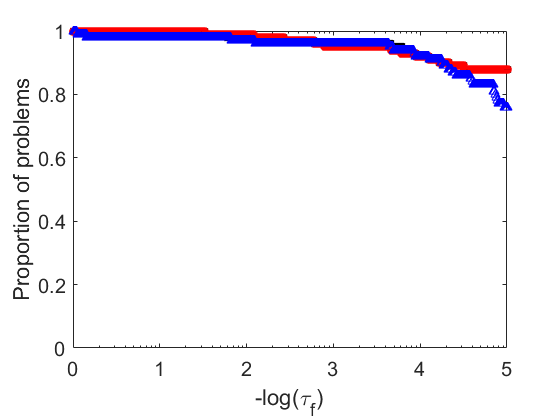}}\\
  \subfloat[][L96, $50\%$]{\includegraphics[width=0.25\textwidth]{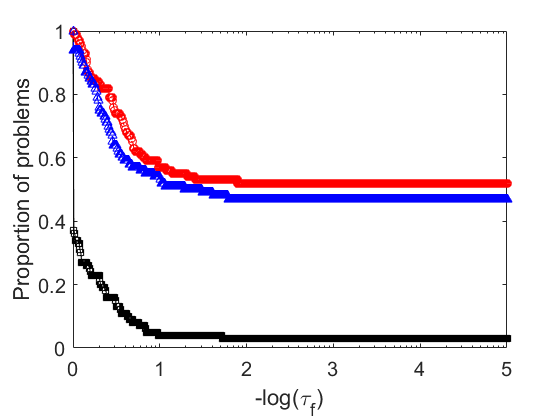}}
 \subfloat[][L96, $25\%$]{\includegraphics[width=0.25\textwidth]{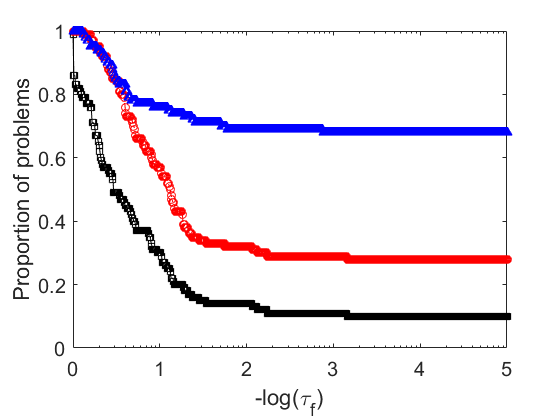}}
 \subfloat[][L96, $10\%$]{\includegraphics[width=0.25\textwidth]{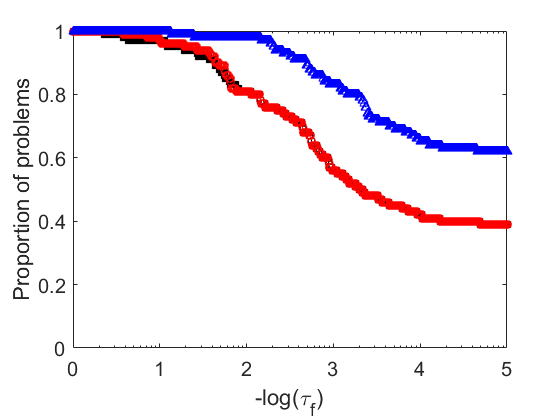}}
 \subfloat[][L96, $5\%$]{\includegraphics[width=0.25\textwidth]{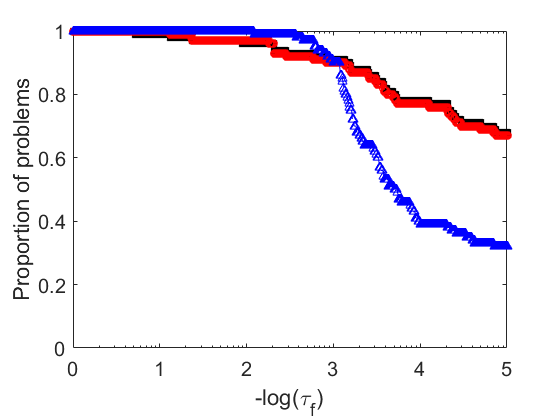}}
 \caption{Accuracy profiles for the GN (black), LS (red) and REG (blue) methods applied to the L63 problems in (a)-(d) and the L96 problems in (e)-(h) where $n_r = 100$, $\tau_e = 8$ and where there is one observation at the end of the time-window. The observation error is $10\%$ and the background error is varied above and below this, as indicated in the plot captions. The GN line is below the LS line in (c), (d), (g) and (h).}
 \label{fig:L63acc}
\end{figure}

\begin{figure}[h!]
\centering
 \subfloat[][L63, $50\%$]{\includegraphics[width=0.25\textwidth]{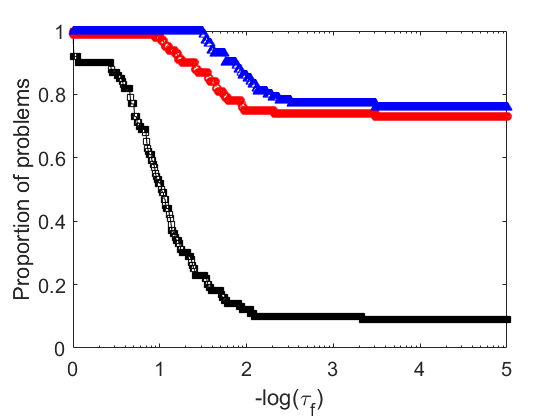}}
 \subfloat[][L63, $25\%$]{\includegraphics[width=0.25\textwidth]{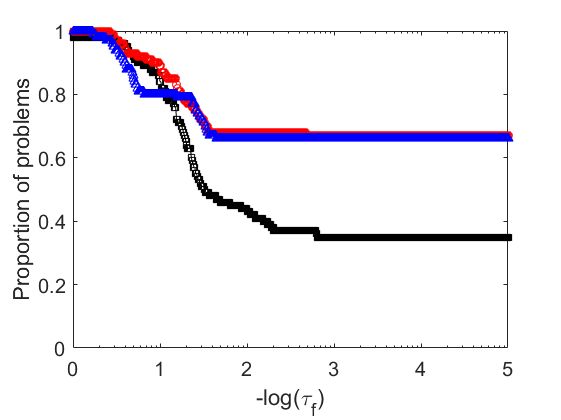}}
 \subfloat[][L63, $10\%$]{\includegraphics[width=0.25\textwidth]{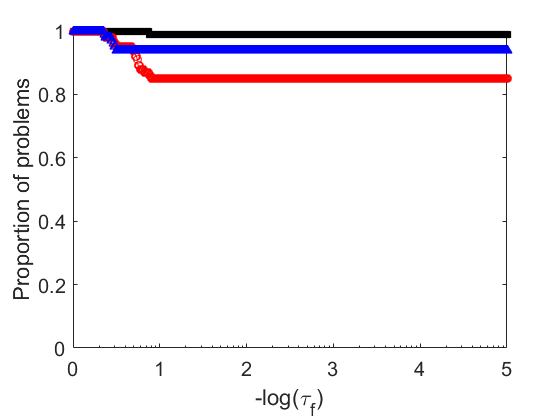}}
 \subfloat[][L63, $5\%$]{\includegraphics[width=0.25\textwidth]{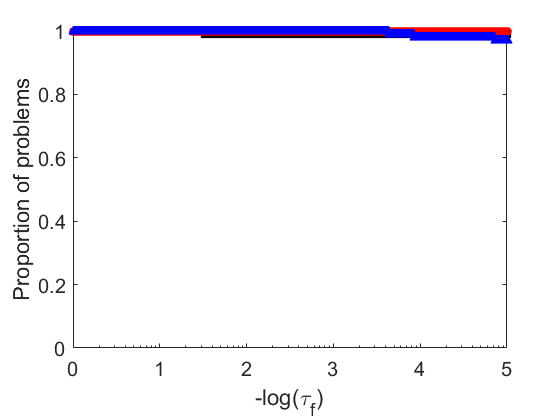}}\\
 \subfloat[][L96, $50\%$]{\includegraphics[width=0.25\textwidth]{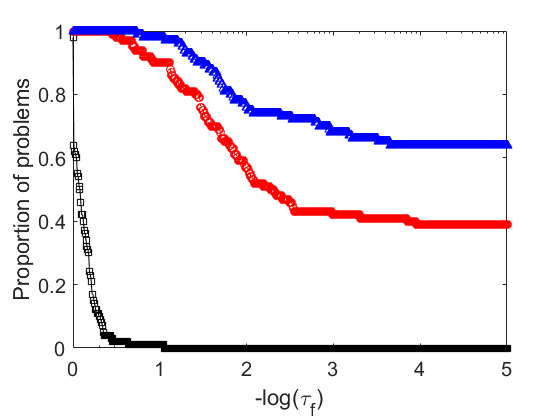}}
 \subfloat[][L96, $25\%$]{\includegraphics[width=0.25\textwidth]{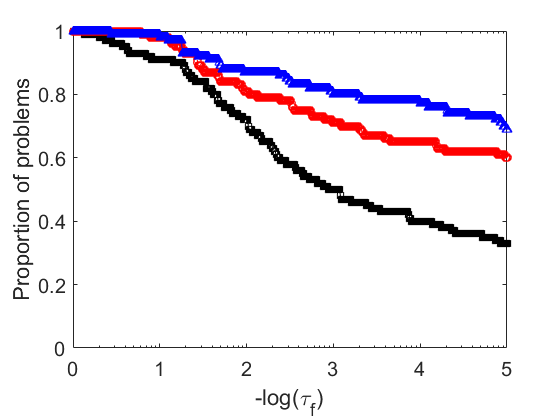}}
 \subfloat[][L96, $10\%$]{\includegraphics[width=0.25\textwidth]{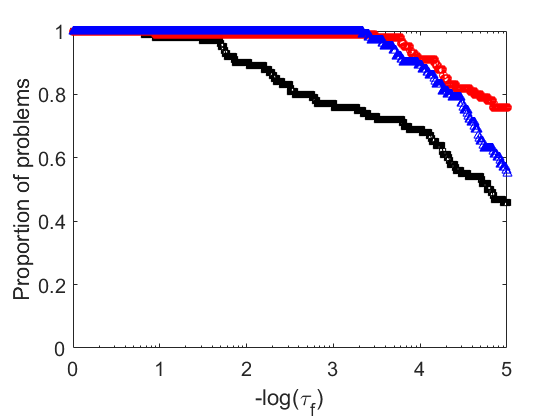}}
 \subfloat[][L96, $5\%$]{\includegraphics[width=0.25\textwidth]{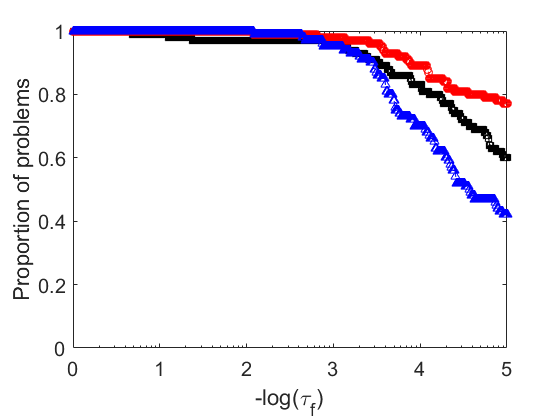}}
 \caption{Accuracy profiles for the GN (black), LS (red) and REG (blue) methods applied to the L63 problems where $\tau_e = 1000$ in (a)-(d) and the L96 problems where $\tau_e = 100$ in (e)-(h). We set $n_r = 100$ and there is one observation at the end of the time-window. The observation error is $10\%$ and the background error is varied above and below this, as indicated in the plot captions.}
 \label{fig:L63acc1000}
\end{figure}

In this section, we study the effect on the performance of the three methods when the uncertainty in the background information is increased whilst the uncertainty in the observations is fixed. Figure \ref{fig:L63acc} shows the accuracy profiles used to benchmark the performance of the GN, LS and REG methods as the tolerance $\tau_f$ is reduced, where $\tau_e = 8$, while Figure \ref{fig:L63acc1000} allows $\tau_e$ to increase for both models with the increase chosen relative to the dimension of the models, i.e. a larger increase in $\tau_e$ is allowed for the L63 problems, where $n=3$, than the L96 problems, where $n=40$. From both these figures, we generally see that as the error in the background is reduced, the performance of all three methods improves. The conditioning of the problem has been shown to depend on the ratio of the standard deviations of the background and observation errors \eqref{sigmaratio} \cite{haben2011conditioning, haben2011conditioninginc}. Therefore, the estimate obtained by any of the optimisation methods may not be accurate enough to produce a reliable forecast if the ratio \eqref{sigmaratio} is large. The accuracy of the estimate obtained by each method will be investigated further later on in the paper.\\ 

Figures \ref{fig:L63acc}(a) and \ref{fig:L63acc}(e) show that a globally convergent method is able to find a smaller function value than GN. As the ratio \eqref{sigmaratio} is reduced, from Figures \ref{fig:L63acc}(b), \ref{fig:L63acc}(c), \ref{fig:L63acc}(f) and \ref{fig:L63acc}(g) we see that the REG method is able to solve the most problems at the highest level of accuracy. When there is less uncertainty in the background versus the observations, Figure \ref{fig:L63acc}(d) shows that for the L63 problems, all three methods are solving close to all of the problems within a high level of accuracy. This is because the three methods are able to solve a large portion of the cases when the problem is well-conditioned, which could explain this result. However, for the L96 problems Figure \ref{fig:L63acc}(h) shows that the GN and LS methods are solving the majority of the problems and REG is not performing as well at higher levels of background accuracy. We can see the performance of REG improving for the L96 problems when more evaluations are allowed in Figure \ref{fig:L63acc1000}(h).\\ 

In Figure \ref{fig:L63acc1000}, where more evaluations are allowed than in Figure \ref{fig:L63acc}, we see a much greater difference between the globally convergent methods and GN when the background error is larger than the observation error. In Figures \ref{fig:L63acc1000}(a), \ref{fig:L63acc1000}(b), \ref{fig:L63acc1000}(e) and \ref{fig:L63acc1000}(f), it appears that when more evaluations are allowed, the performance of GN worsens relative to LS and REG in the case when $\sigma_b$ is large. The globally convergent methods are able to locate estimates of the initial states for the preconditioned 4D-Var problem, which when compared to GN, better minimise the objective function \eqref{4dvarp}. When the background error is the same as the observation error in Figure \ref{fig:L63acc1000}(c), it is GN that is performing better than LS and REG for the L63 problems. For LS, this could be because LS is unnecessarily shortening the GN step, causing slower convergence. For the REG method, the regularisation parameter must be shrunk and therefore, REG requires more iterations to benefit from GN's fast convergence property. \\ 

In Figure \ref{fig:L63acc1000}(d), all three methods are solving essentially the same number of problems, with a slight decrease in success for REG, that again could be due to the need to adjust the regularisation parameter. For the L96 problems, we see a slightly different result. Figures \ref{fig:L63acc1000}(g) and \ref{fig:L63acc1000}(h) show that a globally convergent method is solving more problems, more accurately than GN despite the background error being at most equal to the observation error. This is an interesting result for this higher-dimensional model as we would expect GN to locally converge at a faster rate than the globally convergent methods due to the fact that GN does not need to adjust any parameters; however, we find this not to be the case. \\

In DA, we are interested in knowing the accuracy of the estimate obtained as in applications such as NWP, the estimate is used as the initial conditions for a forecast and so the quality of this forecast will depend on the errors in the estimate. In the following section, we quantify and compare the errors in the estimates obtained by each method.

\subsection{Quality of the analysis}
\label{qualan}

We recall that the initial guess of the algorithms is the reference state $\mathbf{x}_0^{ref}$ perturbed by the background error $\mathbf{\varepsilon_b}$. In order to compare the quality of the estimate obtained by each method, we compare their estimate to the reference state $\mathbf{x}_0^{ref}$ to understand how far the estimates obtained by the methods have deviated from this. The analysis error for each state variable is given by $\varepsilon^a_i = x^a_i - x^{ref}_i$. For each realisation, we calculate the root mean square error (RMSE) of the analysis error, which is the difference between the reference state and the estimate obtained by each method, 
\begin{equation}
    \label{RMSEeq}
    RMSE = \frac{1}{\sqrt{n}} \|\varepsilon^a\|_2.
\end{equation}
For each method, we plot the percentage of problems solved (according to the criterion \eqref{funcxastop} where $\tau_f = 10^{-3}$) within a specified tolerance of the RMSE \eqref{RMSEeq}. We acknowledge in this work that the code for the RMSE profiles has been adapted from the code for the data profiles used in \cite{more2009benchmarking}.\\ 

\begin{figure}[h!]
\centering
 \subfloat[][L63, $50\%$]{\includegraphics[width=0.25\textwidth]{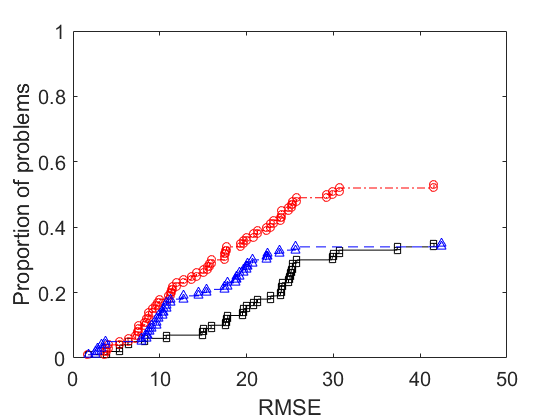}}
 \subfloat[][L63, $25\%$]{\includegraphics[width=0.25\textwidth]{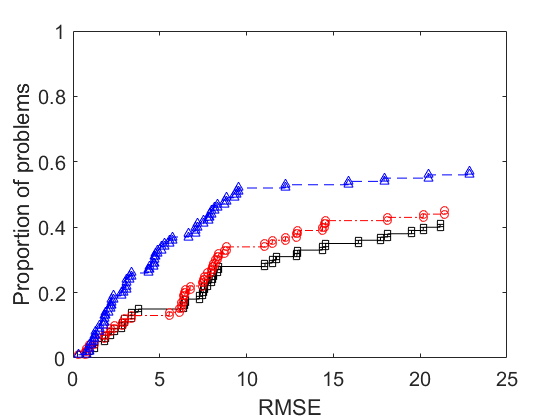}}
 \subfloat[][L63, $10\%$]{\includegraphics[width=0.25\textwidth]{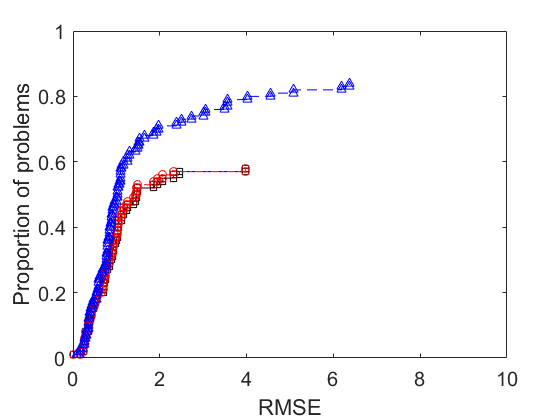}}
 \subfloat[][L63, $5\%$]{\includegraphics[width=0.25\textwidth]{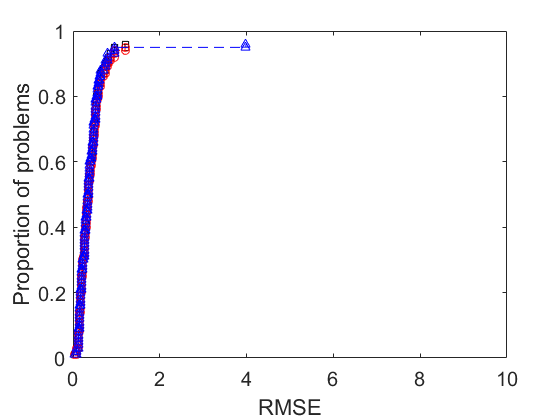}}\\
  \subfloat[][L96, $50\%$]{\includegraphics[width=0.25\textwidth]{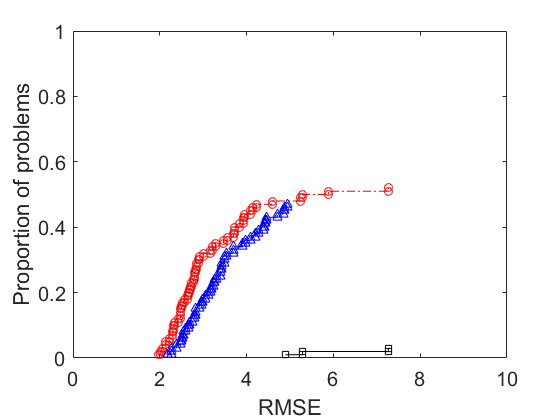}}
 \subfloat[][L96, $25\%$]{\includegraphics[width=0.25\textwidth]{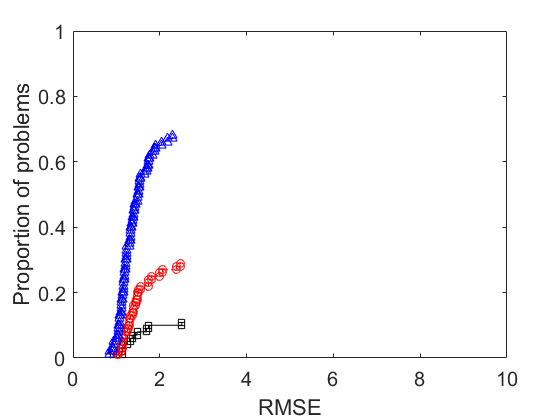}}
 \subfloat[][L96, $10\%$]{\includegraphics[width=0.25\textwidth]{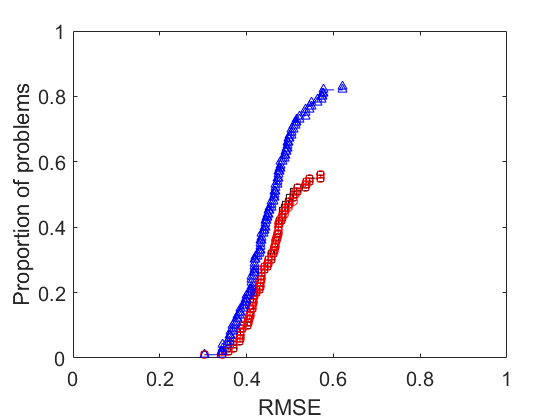}}
 \subfloat[][L96, $5\%$]{\includegraphics[width=0.25\textwidth]{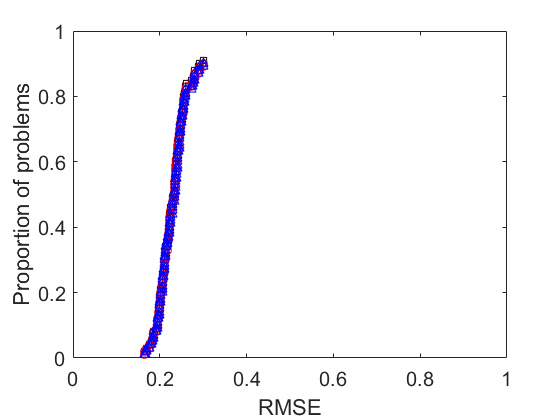}}
  \caption{RMSE plots for the GN (black), LS (red) and REG (blue) methods applied to the L63 problems in (a)-(d) and the L96 problems in (e)-(h) where $n_r = 100$, $\tau_e = 8$, $\tau_f = 10^{-3}$ and where there is one observation at the end of the time-window. The observation error is $10\%$ and the background error is varied above and below this, as indicated in the plot captions.}
  \label{fig:L63RMSE}
\end{figure}

The results for the L63 and L96 problems are in Figure \ref{fig:L63RMSE}, which coincides with the case shown in Figure \ref{fig:L63acc} where $\tau_f = 10^{-3}$. From this, we see that the GN method solves fewer problems within the same level of RMSE accuracy as LS and REG when the background error is large in Figures \ref{fig:L63RMSE}(a), \ref{fig:L63RMSE}(b), \ref{fig:L63RMSE}(e) and \ref{fig:L63RMSE}(f). Furthermore, we see how the RMSE of the analyses successfully found by each method reduces as the background error variance is reduced. This can be seen in the scale of the x axis in Figures \ref{fig:L63RMSE}(a), \ref{fig:L63RMSE}(b), \ref{fig:L63RMSE}(c) and \ref{fig:L63RMSE}(d) for the L63 problems and Figures \ref{fig:L63RMSE}(e), \ref{fig:L63RMSE}(f), \ref{fig:L63RMSE}(g) and \ref{fig:L63RMSE}(h) for the L96 problems. For both models, the concentration of points in Figures \ref{fig:L63RMSE}(a) and \ref{fig:L63RMSE}(e) shows us that the LS method is solving more problems than GN and REG within the same RMSE tolerance. A similar result can be seen for REG in Figures \ref{fig:L63RMSE}(b), \ref{fig:L63RMSE}(c), \ref{fig:L63RMSE}(f) and \ref{fig:L63RMSE}(g). In Figures \ref{fig:L63RMSE}(d) and \ref{fig:L63RMSE}(h), we see that all three methods are performing similarly, the RMSE errors for each of the analyses are very close together. \\

Including more observations constrains the solution to be closer to the reference state when the observation error is small. We next show the effect on the performance of the methods as we include more observations and see if this gives any improvement in the performance of the methods when the background error is much larger than the observation error. 

\subsection{Effect of observations}
\label{effobs}

Within this section, we show how the use of more observation locations in time affects the performance of the three methods. We take the worst case for the three methods when there is a $50\%$ error in the background and see if including more observations in time with a $10\%$ error affects the performance of the methods. For both models, we consider only equally spaced observations in time, one set of observations at time $N$ (Nobs1), times $N/2$ and $N$ (Nobs2), times $N/4, N/2, 3N/4$ and $N$ (Nobs3) and the even time points (Nobs4), where $N=40$. For the Nobs1 case, observations are based on the reference state at the end of the time-window and more observations are included over time in the Nobs2, Nobs3 and Nobs4 cases. This not only increases the condition number of the problem but also constrains the estimate more tightly to the reference state.\\

\begin{figure}[h!]
\centering
 \subfloat[][L63, Nobs1]{\includegraphics[width=0.25\textwidth]{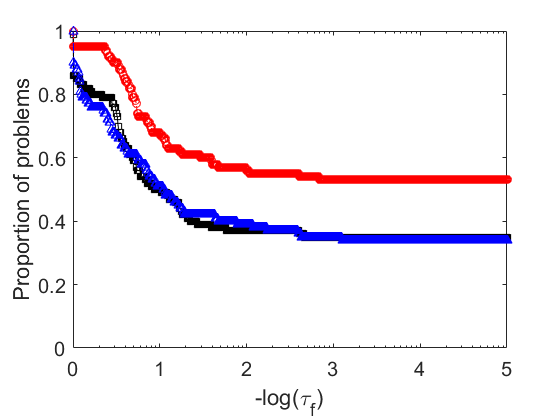}}
 \subfloat[][L63, Nobs2]{\includegraphics[width=0.25\textwidth]{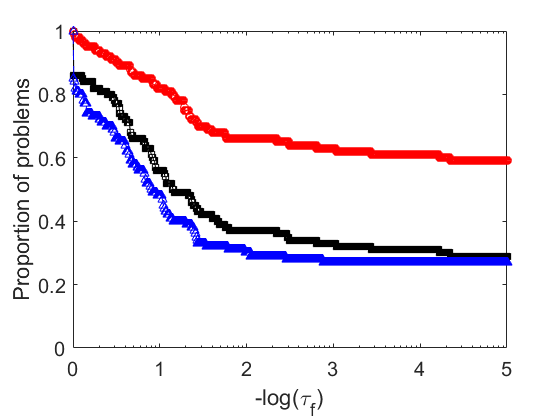}}
 \subfloat[][L63, Nobs3]{\includegraphics[width=0.25\textwidth]{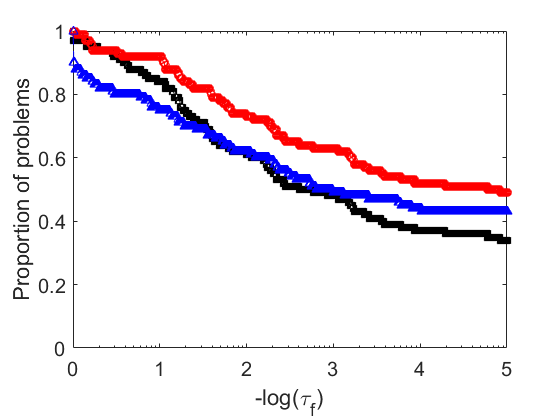}}
 \subfloat[][L63, Nobs4]{\includegraphics[width=0.25\textwidth]{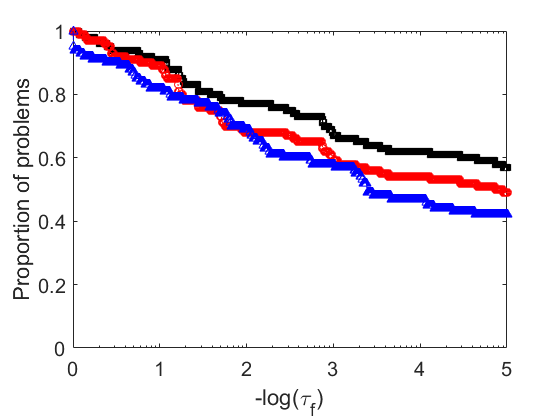}}\\
  \subfloat[][L96, Nobs1]{\includegraphics[width=0.25\textwidth]{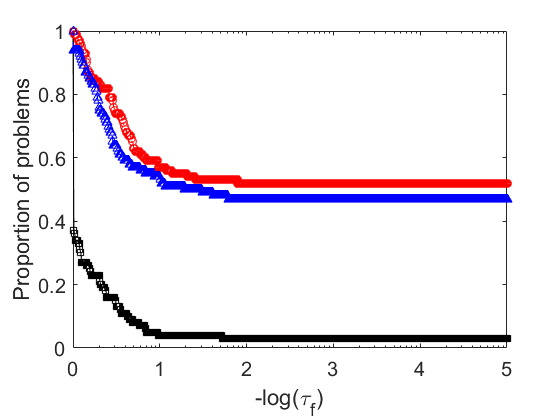}}
 \subfloat[][L96, Nobs2]{\includegraphics[width=0.25\textwidth]{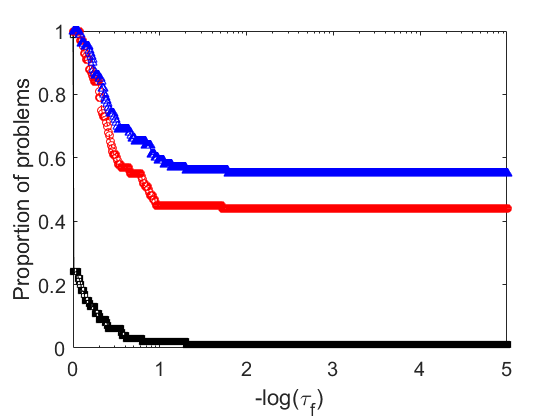}}
 \subfloat[][L96, Nobs3]{\includegraphics[width=0.25\textwidth]{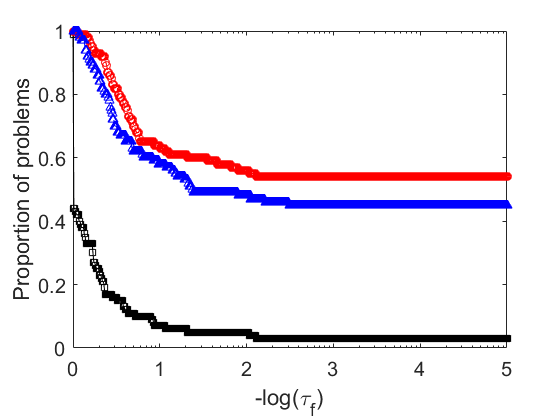}}
 \subfloat[][L96, Nobs4]{\includegraphics[width=0.25\textwidth]{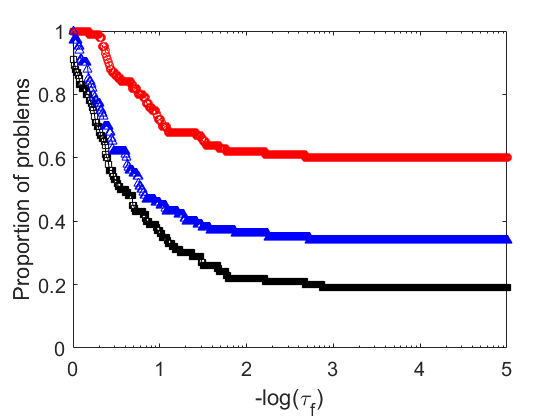}}
 \caption{Accuracy profiles where $n_r = 100$ and $\tau_e = 8$ for the L63 problems in (a)-(d) and the L96 problems in (e)-(h) for different observation locations in time, as indicated in the plot captions, where the background error is $50\%$ and the observation error is $10\%$.}
 \label{fig:nobsL63acc}
\end{figure}

\begin{figure}[h!]
\centering
 \subfloat[][L63, Nobs1]{\includegraphics[width=0.25\textwidth]{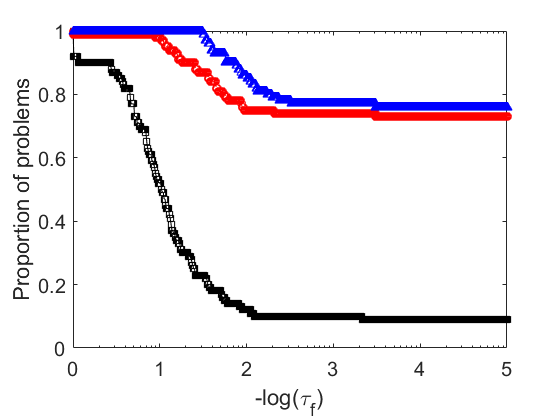}}
 \subfloat[][L63, Nobs2]{\includegraphics[width=0.25\textwidth]{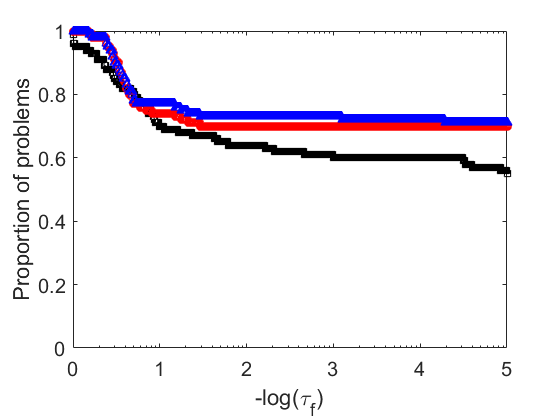}}
 \subfloat[][L63, Nobs3]{\includegraphics[width=0.25\textwidth]{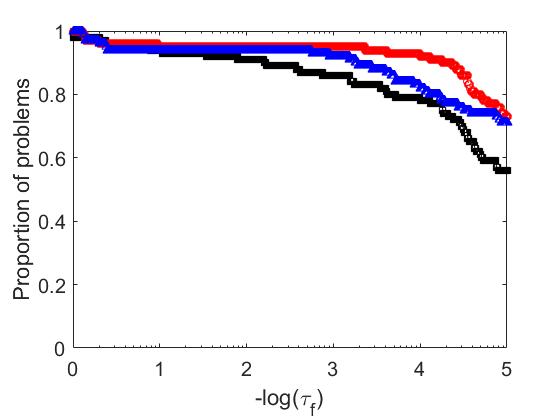}}
 \subfloat[][L63, Nobs4]{\includegraphics[width=0.25\textwidth]{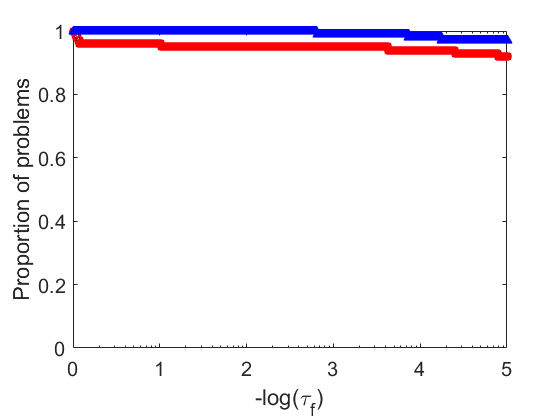}}\\
  \subfloat[][L96, Nobs1]{\includegraphics[width=0.25\textwidth]{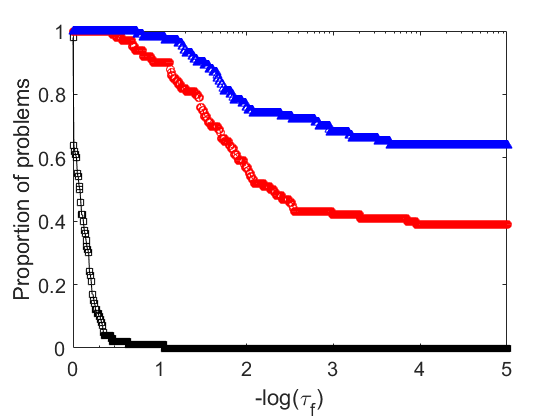}}
 \subfloat[][L96, Nobs2]{\includegraphics[width=0.25\textwidth]{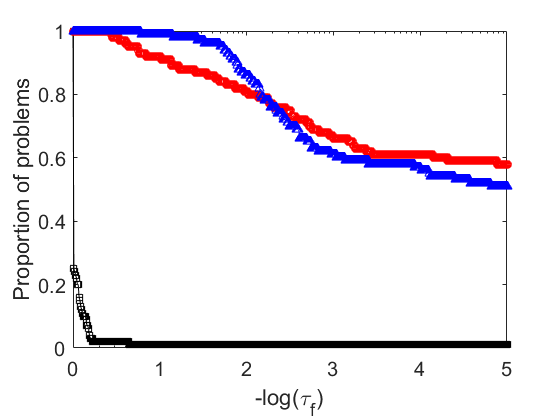}}
 \subfloat[][L96, Nobs3]{\includegraphics[width=0.25\textwidth]{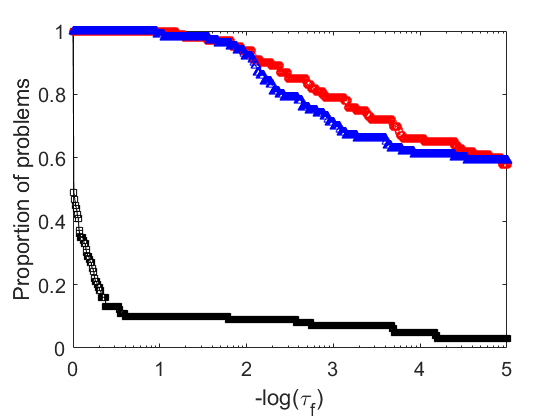}}
 \subfloat[][L96, Nobs4]{\includegraphics[width=0.25\textwidth]{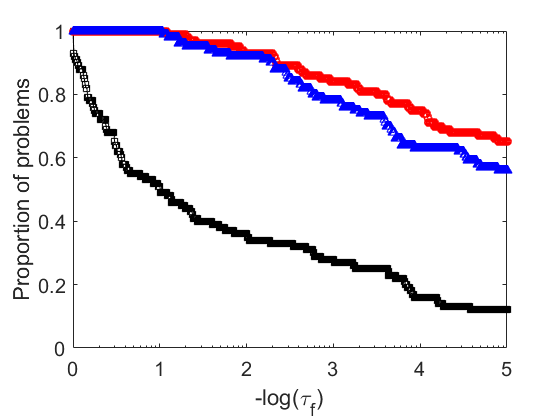}}
 \caption{Accuracy profiles where $n_r = 100$ for the L63 problems where $\tau_e = 1000$ in (a)-(d) and the L96 problems where $\tau_e = 100$ in (e)-(h) for different observation locations in time, as indicated in the plot captions, where the background error is $50\%$ and the observation error is $10\%$. The GN line is below the LS line in (d).}
 \label{fig:nobsL63acc1000}
\end{figure}

For the L63 problems from Figures \ref{fig:nobsL63acc}(a), \ref{fig:nobsL63acc}(b), \ref{fig:nobsL63acc}(c) and \ref{fig:nobsL63acc}(d), we see that as the number of observation locations in time is increased, all three methods are solving more problems at a higher level of accuracy. This is more apparent when more evaluations are allowed as shown in Figure \ref{fig:nobsL63acc1000}(a), \ref{fig:nobsL63acc1000}(b), \ref{fig:nobsL63acc1000}(c) and \ref{fig:nobsL63acc1000}(d). Here, the performance of GN improves drastically between the Nobs1 and Nobs2 cases (Figures \ref{fig:nobsL63acc}(a) and \ref{fig:nobsL63acc}(b)) while there is less significant change in the behaviour of LS and REG. In Figure \ref{fig:nobsL63acc}(d), we see that GN is able to solve more problems than LS and REG. Again, this could be because the LS and REG methods require more iterations to converge when GN is performing well due to the need to adjust their parameters. This argument coincides with Figure \ref{fig:nobsL63acc1000}(d) where more evaluations are allowed and the LS and REG methods are able to perform as well as or better than GN. 
For the L96 problems, we see a different result. From Figure \ref{fig:nobsL63acc}, we only see a significant improvement in the performance of GN in the Nobs4 case (Figure \ref{fig:nobsL63acc}(h)). Otherwise, there is little effect. This conclusion can also be drawn from Figure \ref{fig:nobsL63acc1000}(g) and \ref{fig:nobsL63acc1000}(h) where more evaluations are allowed.\\

Similar studies were carried out on the performance of GN, LS and REG when applied to the preconditioned 4D-Var problem where we instead choose $\mathbf{B} = \sigma_b^2 \mathbf{C}_B$, where $\mathbf{C}_B$ is a correlation matrix; similar conclusions are drawn but due to space constraints, are not included within this paper.

\section{Conclusion} 
\label{concsec}
We have shown that the globally convergent methods, LS and REG, have the capacity to improve current estimates of the DA analysis within the limited time and cost available in DA, through the use of safeguards within GN which guarantee the convergence of the method from any initial guesses. \\

Using the L63 and L96 models in the preconditioned 4D-Var framework, we have shown that when there is more uncertainty in the background information compared to the observations, the GN method may fail to converge in the long time-window case yet the globally convergent methods LS and REG are able to improve the estimate of the initial state. We compare the quality of the estimate obtained using the RMSE of the analysis and show that even in the case where the background is highly inaccurate compared to the observations, the globally convergent methods find estimates with an RMSE less than or equal to the RMSE of the estimates GN obtains. We take the case where the background is highly inaccurate compared to the observations and find that the convergence of all three methods is improved when more observations are included along the time-window. In addition to the numerical results, the assumptions made in the global convergence theorems of both LS and REG when applied to a general nonlinear least-squares problem and a discussion as to whether these assumptions are satisfied in DA is presented in the appendix. We note that preconditioning the second derivative matrix is not necessary for these results to hold, although this is the case we have focused on within our work. \\ 

Our findings are important in DA as they show that in cases where the accuracy of the prior information is poor and when there is limited computational budget, the globally convergent methods are able to minimise the 4D-Var objective function, unlike GN. We recommend that these methods are tested on DA problems with realistic models and for different applications to understand if these conclusions continue to hold. In particular, one should consider such problems where an accurate initial guess for the algorithms is unavailable and a long assimilation time-window is used, as we found that it is in this case that LS and REG have an advantage over GN.\\ 

Within this paper, the 4D-Var inner loop problem is solved exactly. In practice this must be solved inexactly, due to the size of the control vector, and by the use of approximations to meet the computational and time constraints. This is a common area of research in the DA community in order to improve the quality of the assimilation analysis as well as the speed of convergence of the algorithms. Furthermore, in the case where GN performs better than LS and REG, further research is needed on updating the globalisation parameters (stepsize $\alpha^{(k)}$ and regularisation parameter $\gamma^{(k)}$) to speed up convergence. \\

\textbf{Acknowledgements} This work has been funded in part by the UK Engineering and Physical Sciences Research Council Centre for Doctoral Training in Mathematics of Planet Earth, the University of Reading EPSRC studentship (part of Grant/Award Number: EP/N509723/1) and by the NERC National Centre for Earth Observation. We acknowledge in this work that the code for the Lorenz 1996 model was developed by Adam El-Said.\\

\textbf{Declarations of interest} None.

 \section{Appendix}
\begin{appendix}
\section{Convergence theorems}
\label{convthmappend}
In this section, we outline some existing global convergence results for the LS and REG methods and discuss whether the assumptions made hold in DA. We first state the definitions of a local and global minimum of an optimisation problem $\min_{\mathbf{v} \in \mathbb{R}^n} f(\mathbf{v})$ where $f: \mathbb{R}^n \rightarrow \mathbb{R}$ and $\mathbf{v} \in \mathbb{R}^n$.

\begin{definition}[Local minimiser \cite{nocedal2006numerical}] 
A point $\mathbf{v}^*$ is a local minimiser of $f$ if there is a neighbourhood $\mathcal{N}$ of $\mathbf{v}^*$ such that $f(\mathbf{v}^*) \leq f(\mathbf{v})$ for all $\mathbf{v} \in \mathcal{N}$. 
\end{definition}

\begin{definition}[Global minimiser \cite{nocedal2006numerical}]
A point $\mathbf{v}^*$ is a global minimiser of $f: \mathbb{R}^n \rightarrow \mathbb{R}$ if $f(\mathbf{v}^*) \leq f(\mathbf{v})$ for all $\mathbf{v} \in \mathbb{R}^n$.
\end{definition}

A global solution is difficult to locate in most cases due to the nonlinearity of the problems. Therefore, a local solution is often sought by algorithms for nonlinear optimisation. \\

We focus on nonlinear least-squares optimisation problems of the form \eqref{resform} for the remainder of this section. The GN method can only guarantee local convergence under certain conditions and not necessarily global convergence. This is dependent on how close the initial guess is from the local minimum the algorithm locates and whether or not the residual vector $\mathbf{r}$ of \eqref{resform} is a zero vector at a solution $\mathbf{v}^*$. Furthermore, the region of local convergence depends on problem constants not known a priori, such as Lipschitz constants of the gradient. \\

A local convergence result for the GN method can be found in Theorem 10.2.1 of \cite{dennis1996numerical} where the performance of GN is shown to be dependent on whether or not the second-order terms in \eqref{4dvarHess} evaluated at the solution $\mathbf{v}^{*}$ are close to zero. Another local convergence result can be found in Theorem 4 of \cite{gratton2007approximate} where GN is treated as an inexact Newton method. The theorem guarantees convergence of the GN method if for each iteration $k=0, 1, \hdots,$ the norm of the ratio of $\mathbf{Q}(\mathbf{v}^{(k)})$ and $\mathbf{J}(\mathbf{v}^{(k)})^T\mathbf{J}(\mathbf{v}^{(k)})$, the second and first terms of \eqref{4dvarHess} respectively, is less than or equal to some constant $\Hat{\eta}$ where $0 \leq \Hat{\eta} \leq 1$. \\

It is important to note here that the globally convergent methods we are concerned with, namely LS and REG, can only guarantee global convergence to a local minimum under certain conditions and not necessarily to a global minimum. \\

Before we list the assumptions for the global convergence theorems, we first state the definition of the Lipschitz continuity property of a general function $g$ as this is widely used in the theorems.
\begin{definition}[Lipschitz continuous function (see \cite{nocedal2006numerical} A.42)]
Let $g$ be a function where $g : \mathbb{R}^n \rightarrow \mathbb{R}^m$ for general $n$ and $m$. The function $g$ is said to be Lipschitz continuous on some set $\mathcal{N} \subset \mathbb{R}^n$ if there exists a constant $L>0$ such that, 
\begin{equation}
\|g(\mathbf{v}) - g(\mathbf{w})\| \leq L\|\mathbf{v} - \mathbf{w}\|, \qquad \forall \mathbf{v},\mathbf{w} \in \mathcal{N}.
\end{equation}
\end{definition}

The following assumptions are used to prove global convergence of both the LS and REG methods. 

\begin{assump}
\label{rboundedaboveassump}
$\mathbf{r}$ is uniformly bounded above by $\omega >0$ such that $\|\mathbf{r}(\mathbf{v})\| \leq \omega$.
\end{assump}

\begin{assump}
\label{rlipcontassump}
$\mathbf{r} \in \mathcal{C}^1(\mathbb{R}^n)$ is Lipschitz continuous on $\mathbb{R}^n$ with Lipschitz constant $L_r > 0$.
\end{assump}

\begin{assump}
\label{Jaclipcontassump}
$\mathbf{J}$ is Lipschitz continuous on $\mathbb{R}^n$ with Lipschitz constant $L_J > 0$.
\end{assump}

We remark that for the LS method, we can weaken assumptions \ref{rlipcontassump} and \ref{Jaclipcontassump} using the open set $\mathcal{N}$ containing the level set 
\begin{equation}
\label{levelset}
\mathcal{L} = \left\lbrace \mathbf{v} \in \mathbb{R}^n| \mathcal{J}(\mathbf{v}) \leq \mathcal{J}(\mathbf{v}^{(0)})\right\rbrace.
\end{equation}

In order to achieve the sufficient decrease property of the LS method, the following assumption must be made.
\begin{assump}
\label{JacfullrankLSassump}
$\mathbf{J}(\mathbf{v})$ in \eqref{precondresJac} is uniformly full rank for all $\mathbf{v} \in \mathbb{R}^n$, that is, the singular values of $\mathbf{J}(\mathbf{v})$ are uniformly bounded away from zero, so there exists a constant $\nu$ such that $\| \mathbf{J}(\mathbf{v}) \mathbf{z} \|  \geq \nu \| \mathbf{z} \|$ for all $\mathbf{v}$ in a neighbourhood $\mathcal{N}$ of the level set $\mathcal{L}$ where $\mathbf{z} \in \mathbb{R}^n$.
\end{assump}

In 4D-Var practice, it is reasonable to assume that the physical quantities are bounded. Therefore, we can say that both $\mathbf{x}_0 - \mathbf{x}^b$ and the innovation vector $\mathbf{y} - \mathcal{H}(\mathbf{x})$ are bounded in practice, thus satisfying assumption \ref{rboundedaboveassump}.
In 4D-Var, we must assume that the nonlinear model $\mathcal{M}_{0,i}$ is Lipschitz continuous in order for \ref{rlipcontassump} to hold. As discussed in \cite{moodey2013nonlinear}, this is a common assumption made in the meteorological applications. However, we cannot say that this is necessarily the case in 4D-Var practice.\\

In order for the Jacobian $\mathbf{J}$ to be Lipschitz continuous, we require its derivative to be bounded above by its Lipschitz constant. Therefore, for assumption \ref{Jaclipcontassump} to hold, we require $\mathbf{r}$ to be twice continuously differentiable in practice, which is a common assumption made in 4D-Var, and also, that these derivatives of $\mathbf{r}$ are bounded above. \\

As mentioned in Section \ref{vardasection}, the preconditioned 4D-Var Hessian \eqref{precondhess} is full rank by construction as it consists of the identity matrix and a non-negative definite term. Therefore, the Jacobian of the residual of the preconditioned problem in \eqref{precondresJac} is full rank and assumption \ref{JacfullrankLSassump} holds. This is also the case for the standard 4D-Var problem \eqref{4dvar}, because of the presence of $\mathbf{B}^{1/2}$ in its Jacobian. \\

We now outline the global convergence theorems for the LS and REG methods, using these assumptions.

\subsection{Global convergence of the LS method}
\label{LSappendix}

Nocedal et al. outline the proof for the GN method with Wolfe line search conditions in \cite{nocedal2006numerical}, which uses the Zoutendijk condition. This proof can be adapted to prove the global convergence theorem of the LS method, Algorithm \ref{LSalg}, given as follows. 

\begin{theorem}[Global convergence for the Gauss-Newton with bArmijo line search method, Algorithm \ref{LSalg}]
\label{GNarmijotheorem}
Suppose we have a function $\mathcal{J} = \frac{1}{2}\mathbf{r}^T\mathbf{r}$ and its gradient $\nabla \mathcal{J} = \mathbf{J}^T\mathbf{r}$ where $\mathbf{r} \in \mathcal{C}^1(\mathbb{R}^n)$ and $\mathbf{J}$ is the Jacobian of $\mathbf{r}$. Assume \ref{rboundedaboveassump} - \ref{JacfullrankLSassump} hold. Then if the iterates $\lbrace\mathbf{v}^{(k)}\rbrace$ are generated by the GN method with stepsizes $\alpha^{(k)}$ that satisfy the Armijo condition \eqref{bArmijoalg}, we have
\begin{equation}
\label{LSthmeq}
\lim_{k \rightarrow \infty} \mathbf{J}(\mathbf{v}^{(k)})^T \mathbf{r}(\mathbf{v}^{(k)}) = 0.
\end{equation}
That is, the gradient norms converge to zero, and so the Gauss-Newton method with bArmijo line search is globally convergent.
\end{theorem}

The proof of Theorem \ref{GNarmijotheorem} requires the bArmijo chosen stepsizes $\alpha^{(k)}$ to be bounded below, which can be derived using assumptions \ref{rboundedaboveassump} - \ref{Jaclipcontassump}. Using this lower bound, as well as assumption \ref{JacfullrankLSassump}, we are able to prove the Zoutendijk condition (as in \cite{nocedal2006numerical}) and its variant 
\begin{equation}
\label{bArmijolbcase2res}
\sum_{k \geq 0} \cos(\theta^{(k)})\|\nabla \mathcal{J}(\mathbf{v}^{(k)})\|_2\|\mathbf{s}^{(k)}\|_2 < \infty
\end{equation}
hold. Both the Zoutendijk condition and its variant \eqref{bArmijolbcase2res} use the angle between $\mathbf{s}^{(k)}$ (the GN search direction) and $-\nabla \mathcal{J}(\mathbf{v}^{(k)})$ (the steepest descent direction), $\theta^{(k)}$, which is given by
\begin{equation}
\label{LSangle}
\cos(\theta^{(k)}) = \frac{(-\nabla \mathcal{J}(\mathbf{v}^{(k)}))^T\mathbf{s}^{(k)}}{\|\nabla \mathcal{J}(\mathbf{v}^{(k)})\|_2 \|\mathbf{s}^{(k)}\|_2}.
\end{equation} 
By showing that the angle is uniformly bounded away from zero with $k$, one can show that GN with line search is a globally convergent method.\\

We will next present the global convergence theorem for the REG method. The REG method has no sufficient decrease condition as in the LS method. Therefore, the use of the level set \eqref{levelset} is not required. The assumptions for convergence are similar to the LS method aside from the requirement of $\mathbf{J}(\mathbf{v})$ being full rank.

\subsection{Global convergence of the REG method}
\label{REGappendix}

The global convergence theorem for the GN with quadratic regularisation method, Algorithm \ref{REGalg}, is given as follows.

\begin{theorem}[Global convergence for the Gauss-Newton with regularisation method, Algorithm \ref{REGalg}]
\label{regglobalproof}
Suppose we have a function $\mathcal{J} = \frac{1}{2}\mathbf{r}^T\mathbf{r}$ and its gradient $\nabla \mathcal{J} = \mathbf{J}^T\mathbf{r}$ where $\mathbf{r} \in \mathcal{C}^1(\mathbb{R}^n)$ and $\mathbf{J}$ is the Jacobian of $\mathbf{r}$. Assume \ref{rboundedaboveassump} - \ref{Jaclipcontassump} hold. Then if the iterates $\lbrace\mathbf{v}^{(k)}\rbrace$ are generated by the Gauss-Newton with regularisation method, we have that
\begin{equation}
\label{REGthmeq}
\lim_{k \rightarrow \infty} \mathbf{J}(\mathbf{v}^{(k)})^T \mathbf{r}(\mathbf{v}^{(k)}) = 0.
\end{equation}
That is, the gradient norms converge to zero, and so the Gauss-Newton method with regularisation is globally convergent.
\end{theorem}

We first note that some adaptations of the lemmas from the global convergence proof of the Adaptive Regularisation algorithm using Cubics (ARC method) are used to prove Theorem \ref{regglobalproof}, see \cite{cartis2011adaptive1} and \cite{cartis2011adaptive2}. We begin the proof by deriving an expression for the predicted model decrease in terms of the gradient. We require the use of an upper bound on $\gamma^{(k)}$, denoted as $\gamma_{\max}$, which is derived using a property of Lipschitz continuous gradients. We show that $\gamma^{(k)} \leq \gamma_{\max}$ for all $k \geq 0$ by first showing that if $\gamma^{(k)}$ is large enough, then we have a successful step so that $\gamma^{(k)}$ can stop increasing due to unsuccessful steps in Algorithm \ref{REGalg}. We use the expression for $\gamma_{\max}$ to prove global convergence of the REG method under assumptions \ref{rboundedaboveassump}-\ref{Jaclipcontassump} by showing that the gradient norms converge to zero as we iterate.\\

Note that for both the LS and REG, if $\mathbf{r}(\mathbf{v}^{(k)}) \rightarrow 0$, i.e. \eqref{resform} is a zero residual problem, then we have that \eqref{LSthmeq} and \eqref{REGthmeq} hold as $| \mathcal{J}(\mathbf{v}^{(k)}) |$ is uniformly bounded. However, in practice the variational problem is not usually a zero residual problem.

\end{appendix}

\bibliography{Bibfile}
\end{document}